\RequirePackage{fix-cm}
\documentclass[smallextended,envcountsect]{svjour_arxiv}       % onecolumn (second format)
\smartqed  % flush right qed marks, e.g. at end of proof
\usepackage[dvips]{graphicx}

\usepackage{here}
\usepackage{amsmath,amssymb,xypic,url}
\usepackage[subrefformat=parens]{subcaption}
\newcommand{\R}{{\mathbb R}}
\newcommand{\N}{{\mathbb N}}
\newcommand{\Hb}{{H^1(\Omega)}}
\newcommand{\Ho}{{H^1_0(\Omega)}}
\newcommand{\HD}{{H^1_{0,D}(\Omega)}}
\newcommand{\Hi}{{H^{-1}(\Omega)}}
\newcommand{\Lo}{{L^2(\Omega)}}
\newcommand{\dvg}{{\operatorname{div}}}
\newcommand{\into}{\int_\Omega}
\newcommand{\intg}{\int_\Gamma}
\newcommand{\eps}{\varepsilon}
\newcommand{\DS}{\displaystyle}
\newcommand{\St}{{\mbox{\tiny$S$}}}
\newcommand{\PP}{{\mbox{\tiny$P\!P$}}}

\spnewtheorem{Th}{Theorem}[section]{\bf}{\it}
\spnewtheorem{Lem}[Th]{Lemma}{\bf}{\it}
\spnewtheorem{Prop}[Th]{Proposition}{\bf}{\it}
\spnewtheorem{Cor}[Th]{Corollary}{\bf}{\it}
\spnewtheorem{Def}[Th]{Definition}{\bf}{\it}
\spnewtheorem{Prob}[Th]{Problem}{\bf}{\it}
\spnewtheorem{Rem}[Th]{Remark}{\bf}{\it}

\renewcommand{\theequation}{\arabic{section}.\arabic{equation}}

\begin{document}

\title{Analysis of a projection method for the Stokes problem
using an $\eps$-Stokes approach
}

\titlerunning{$\varepsilon$-Stokes approach}        % if too long for running head

\author{Masato Kimura \and Kazunori Matsui \and Adrian Muntean\and Hirofumi Notsu}


\institute{
Masato Kimura\at
Faculty of Mathematics and Physics, Kanazawa University, Kanazawa
920-1192, Japan. %\\
%\email{mkimura@se.kanazawa-u.ac.jp}           %  \\
\and 
Kazunori Matsui (Corresponding author)\at
Division of Mathematical and Physical Sciences, \\
Graduate School of Natural Science and Technology, Kanazawa University, \\
Kanazawa 920-1192, Japan, \\
\email{first-lucky@stu.kanazawa-u.ac.jp}
\and 
Adrian Muntean \at 
Department of Mathematics and Computer Science, Karlstad University, \\
Universitetsgatan 2, 651 88 Karlstad Sweden. %\\
%\email{adrian.muntean@kau.se}
\and 
Hirofumi Notsu \at
Faculty of Mathematics and Physics, Kanazawa University, Kanazawa
920-1192, Japan, \\
Japan Science and Technology Agency, PRESTO, Kawaguchi 332-0012, Japan. %\\
%\email{notsu@se.kanazawa-u.ac.jp}
}


\maketitle

\begin{abstract}

We generalize pressure boundary conditions of an $\eps$-Stokes problem.
Our $\eps$-Stokes problem connects the classical Stokes problem 
and the corresponding pressure-Poisson equation using one parameter $\eps>0$.
For the Dirichlet boundary condition, it is proven in K. Matsui and A. Muntean (2018) that
the solution for the $\eps$-Stokes problem converges to the one for the Stokes problem
as $\eps$ tends to 0, and to the one for the pressure-Poisson problem
as $\eps$ tends to $\infty$. Here, we extend these results 
to the Neumann and mixed boundary conditions.
We also establish error estimates in suitable norms
between the solutions to the $\eps$-Stokes problem,
the pressure-Poisson problem and the Stokes problem, respectively.
Several numerical examples are provided to show that 
several such error estimates are optimal in $\eps$.
Our error estimates are improved if one uses the Neumann boundary conditions.
In addition, we show that the solution to the $\eps$-Stokes problem 
has a nice asymptotic structure.

\keywords{Stokes problem \and Pressure-Poisson equation \and
Asymptotic analysis \and Finite element method}
\subclass{76D03 \and 35Q35 \and 35B40 \and 65N30 }
\end{abstract}

%
% ----------------------------------------------------------------------------
\section{Introduction}
\label{intro}
% ----------------------------------------------------------------------------
%

Let $\Omega$ be a bounded domain in $\R^n(n\ge 2,n\in\N)$
with Lipschitz continuous boundary $\Gamma$ and
let $F:\Omega\rightarrow\R^n$ be a given applied force field and 
$u_b:\Gamma\rightarrow\R^n$ be a given Dirichlet boundary data satisfying 
$\int_\Gamma u_b\cdot \nu=0$,
where $\nu$ is the unit outward normal vector on $\Gamma$.
A strong form of the Stokes problem is given as follows. 
Find $u_\St:\Omega\rightarrow\R^n$ and $p_\St:\Omega\rightarrow\R$ such that
\begin{align}\tag{S}
\left\{\begin{array}{ll}
-\Delta u_\St+\nabla p_\St = F & \mbox{in }\Omega, \\
\dvg u_\St = 0 & \mbox{in }\Omega, \\
u_\St=u_b & \mbox{on }\Gamma,
\end{array}\right.
\end{align}
where $u_\St$ and $p_\St$ are the velocity and the pressure of the flow
governed by (S), respectively.
We refer to \cite{Temam} for the details on the Stokes problem
(i.e., physical background and corresponding mathematical analysis).
Taking the divergence of the first equation, we obtain
\begin{align}\label{ppstrong}
\dvg F=\dvg(-\Delta u_\St+\nabla p_\St)=-\Delta (\dvg u_\St)+\Delta p_\St=\Delta p_\St.
\end{align}
This equation is often called the pressure-Poisson equation and 
is used in numerical schemes such as MAC (marker and cell), SMAC 
(simplified MAC) and the projection methods 
(see, e.g., \cite{Amsden_Harlow,Chorin68,Cummins_Rudman,Guermond,mac1,Kim_Moin,mac2,Perot}).
Based on the above, we consider a similar problem.
Find $u_\PP:\Omega\rightarrow\R^n$ and $p_\PP:\Omega\rightarrow\R$ satisfying
\begin{align}\tag{PP}
\left\{\begin{array}{ll}
-\Delta u_\PP+\nabla p_\PP = F & \mbox{in }\Omega, \\
-\Delta p_\PP = -\dvg F& \mbox{in }\Omega, \\
u_\PP=u_b & \mbox{on }\Gamma, \\
+\mbox{boundary condition for }p_\PP. &
\end{array}\right.
\end{align}
We call this problem the pressure-Poisson problem.
The idea of using (\ref{ppstrong}) instead of $\dvg u_\St=0$ is useful
for calculating the pressure numerically in the Navier--Stokes equation.
For example, this idea is used in MAC, SMAC and projection methods.
The Dirichlet boundary condition for the pressure is used in 
an outflow boundary \cite{Chan_Street,Viecelli}. 
See also \cite{Conca_etc94,Conca_etc95,Marusic}.

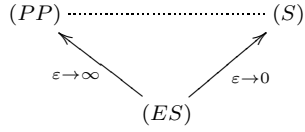
\begin{figure}[H]
\[
\xymatrix{
  (PP)\ar@{.}[rr]& &(S)\\
 &(ES)\ar[ul]^{ \eps\rightarrow\infty}\ar[ur]_{\eps\rightarrow 0}&
}
\]
\caption{
  Sketch of the connections between the problems (S), (PP) and (ES).
}
\label{diagram}
\end{figure}

We introduce an ``interpolation'' between problems (S) and (PP).
For $\eps>0$, find $u_\eps:\Omega\rightarrow\R^n$ and
$p_\eps:\Omega\rightarrow\R$ such that
\begin{align}\tag{ES}
\left\{\begin{array}{ll}
-\Delta u_\eps+\nabla p_\eps = F & \mbox{in }\Omega, \\
-\eps\Delta p_\eps+\dvg u_\eps = -\eps\dvg F& \mbox{in }\Omega, \\
u_\eps=u_b & \mbox{on }\Gamma, \\
+\mbox{boundary condition for }p_\eps. &
\end{array}\right.
\end{align}
This problem is called the $\eps$-Stokes problem (ES) in \cite{prev}.
In \cite{Douglas_Wang,Glowinski,Hughes}, the authors treat this problem
as an approximation of the Stokes problem to avoid numerical instabilities.
The $\eps$-Stokes problem approximates the Stokes problem (S)
as $\eps\rightarrow0$ and the pressure-Poisson problem (PP)
as $\eps\rightarrow\infty$ (Fig.~\ref{diagram}).
It is shown in \cite{prev} that
if $p_\St\in\Hb$ then there exists a constant $c>0$ independent of $\eps$ such that
\[\begin{array}{rrl}
  \|u_\St-u_\PP\|_{\Hb^n}
  &\le& c\|\gamma_0 p_\St-\gamma_0 p_\PP\|_{H^{1/2}(\Gamma)},\\
  \|u_\St-u_\eps\|_{\Hb^n}
  &\le& c\|\gamma_0 p_\St-\gamma_0 p_\PP\|_{H^{1/2}(\Gamma)},
\end{array}\]
where $\gamma_0\in B(\Hb,H^{1/2}(\Gamma))$ is the standard trace operator \cite{Girault}. 
From the first inequality, if we have a good prediction value for pressure on $\Gamma$,
then $u_\PP$ is a good approximation of $u_\St$.
Moreover, $u_\eps$ is also a good approximation of $u_\St$ from the second inequality. 

Next we specify the boundary conditions for $p_\PP$ and $p_\eps$.
We assume that the boundary $\Gamma$ is a union of two open subsets
$\Gamma_D$ and $\Gamma_N$ such that
\[
  |\Gamma\setminus(\Gamma_D\cup\Gamma_N)|=0,\quad
  |\Gamma_D|> 0,\quad
  \Gamma_D\cap\Gamma_N=\emptyset,
\]
and number of connected components of $\Gamma_D$ and $\Gamma_N$
with respect to the relative topology of $\Gamma$ are finite.
We consider a Neumann boundary condition (\ref{bcn}) and
a mixed boundary condition (\ref{bcm}),
\begin{align}\label{bcn}
{\DS \frac{\partial p_\PP}{\partial \nu}=g_b}\mbox{ on }\Gamma,\quad
{\DS \frac{\partial p_\eps}{\partial \nu}=g_b}\mbox{ on }\Gamma,
\end{align}
\begin{align}\label{bcm}
\left\{\begin{array}{ll}
p_\PP=p_b & \mbox{on }\Gamma_D,\\[8pt]
{\DS \frac{\partial p_\PP}{\partial \nu}=g_b} & \mbox{on }\Gamma_N,
\end{array}\right.\quad
\left\{\begin{array}{ll}
  p_\eps=p_b & \mbox{on }\Gamma_D,\\[8pt]
  {\DS \frac{\partial p_\eps}{\partial \nu}=g_b} & \mbox{on }\Gamma_N,
  \end{array}\right.
\end{align}
where $p_b:\Gamma_D\rightarrow\R$ and $g_b=\Gamma\rightarrow\R$ 
satisfying $\intg g_b=\intg\dvg F$ are given boundary data.

In \cite{prev}, the authors impose Dirichlet boundary conditions
for $p_\PP$ and $p_\eps$ 
(i.e., (\ref{bcm}) with $\Gamma_D=\Gamma$ and $\Gamma_N=\emptyset$.) 
For such boundary conditions, they introduce a weak solution 
$(u_\eps,p_\eps)$ to the $\eps$-Stokes problem (ES) and
prove that $(u_\eps,p_\eps)$ strongly converges in $\Hb^n\times\Hb$ to
a weak solution to the pressure-Poisson problem (PP)
as $\eps\rightarrow\infty$ and weakly converges in $\Ho^n\times(\Lo/\R)$
to a weak solution $(u_\St,p_\St)$ to the Stokes problem (S)
as $\eps\rightarrow0$.
Moreover, if $p_\St\in\Hb$, then strong convergence of $(u_\eps,p_\eps)$ 
to $(u_\St,p_\St)$ as $\eps\rightarrow0$ holds.

In this paper, we generalize the Dirichlet boundary condition 
of $p_\PP$ and $p_\eps$ to both the Neumann boundary condition (\ref{bcn}) 
and the mixed boundary condition (\ref{bcm}),
and prove the corresponding convergence result (see Theorem 3.1, 4.2 and 4.3).
Since the mixed boundary condition for pressure often appears 
in engineering problems, this generalization of the boundary conditions
for pressure is important.
In addition, for the Neumann boundary condition,
we estimate the error between the weak solutions to (ES) and (S) provided $p_\St\in\Hb$.
We also give an asymptotic expansion for the weak solution to (ES).
We furthermore check this convergence result using numerical computations.

The organization of this paper is as follows.
In Section \ref{sec_wellposed} we introduce the notation used in this work
and the weak form of these problems. We also 
prove the well-posedness of the problems (PP) and (ES)
and show some their properties.
In Section \ref{sec_ep} we study that the solution to (ES) converges
to the solution to (PP) in the strong topology as $\eps\rightarrow\infty$.
We also explore here the structure of the regular perturbation asymptotics.
Section \ref{sec_es} is devoted to proving that the solution to (ES)
converges to the solution to (S) in the weak and strong topology
as $\eps\rightarrow 0$.
Finally, in Section \ref{sec_numex},
we show several numerical examples of these problems.
The numerical errors between the problems (ES) and (PP),
and between the problems (ES) and (S)
using the P2/P1 finite element method.
Proofs for several theorems which are similar to ones in \cite{prev}
are described.

%
% ----------------------------------------------------------------------------
\section{Well-posedness}\label{sec_wellposed}
% ----------------------------------------------------------------------------
%

In this section, we introduce the notation and the weak form 
of the problems (S), (PP) and (ES), and prove their well-posedness.
We give estimates between these solutions by
using a pressure error on the boundary $\Gamma$.

%
% ----------------------------------------------
\subsection{Notation}\label{notations}
% ----------------------------------------------
%
We set
\[\begin{array}{ll}
C^\infty_0(\Omega)^n
&{\DS :=\{\varphi\in C^\infty(\Omega)^n\,|\,\mbox{supp}(\varphi)\subset\Omega\},}\\[8pt]
\Lo/\R&:={\DS \left\{\psi\in\Lo\,\middle|\,\into \psi=0\right\},}\\[8pt]
\Hb/\R&:=\Hb\cap(\Lo/\R),\\[8pt]
\HD&:={\DS \{\psi\in\Hb\mid\psi|_{\Gamma_D}=0\}.}\\[8pt]
\end{array}\]
For $m\in\N$, $\Hi^m:=(\Ho^m)^*$ is equipped with the dual norm
\[
\|f\|_{\Hi^m}:=\sup_{\varphi\in S_m}\langle f,\varphi\rangle
\mbox{ for }f\in\Hi^m,
\]
where
\[
S_m:=\{\varphi\in\Ho^m~|~\|\nabla\varphi\|_{\Lo^{n\times m}}=1\}.
\]

Let $Q\subset\Hb$ be a closed subspace such that 
there exists a constant $c>0$ for which
$\|q\|_\Lo\le c\|\nabla q\|_{\Lo^n}$ for all $q\in Q$.
The dual space $Q^*$ is equipped with the norm
\[
\|f\|_{Q^*}:=\sup_{\psi\in S_Q}\langle f,\psi\rangle
\]
for $f\in Q^*$,
where 
\[
S_Q:=\{\psi\in Q~|~\|\nabla\psi\|_{\Lo^{n}}=1\}.
\]
We define 
\[\begin{array}{rl}
[p]&:={\DS p-\frac{1}{|\Omega|}\into p,}\\[8pt]
\|p\|_{\Lo/\R}
&:={\DS \inf_{a\in\R}\|p-a\|_\Lo=\|[p]\|_\Lo,}\\[8pt]
\langle\nabla p,\varphi\rangle
&:={\DS-\into p\,\dvg\varphi}\\[8pt]
{\DS \nabla u:\nabla\varphi}
&:={\DS \sum^n_{i=1}\frac{\partial u}{\partial x_i}\cdot\frac{\partial\varphi}{\partial x_i}}
\end{array}\]
for all $p\in\Lo,u,\varphi\in\Ho^n$, 
where $|\Omega|$ is the volume of $\Omega$.

%
% ----------------------------------------------
\subsection{Preliminary results}\label{pre_results}
% ----------------------------------------------
%

Let $\gamma_0\in B(\Hb, H^{1/2}(\Gamma))$ be
the standard trace operator.
The trace operator $\gamma_0$ is surjective and satisfies 
$\mbox{Ker}(\gamma_0)=\Ho$ \cite[Theorem 1.5]{Girault}.
Let $\nu$ be the unit outward normal for $\Gamma$.
Since $\nu$ is a unit vector, 
$\Hb^3\ni u\mapsto u\cdot\nu:=(\gamma_0 u)\cdot\nu\in H^{1/2}(\Gamma)$
is a linear continuous map.
For all $u\in \Hb^3$ and $\omega\in\Hb$, the following Gauss divergence formula holds:
\[
\into u\cdot\nabla \omega+\into(\dvg u)\omega=\int_\Gamma(u\cdot\nu)\omega.
\]

We recall the following four embedding theorems
which plays an important role in the proof of
the existence of pressure solutions to the Stokes problem.
For the proof of Theorem \ref{lem_Necas}, 
see \cite[Lemma 7.1]{Necas} and
\cite[Theorem 3.2 and Remark 3.1]{Duvaut}.
\begin{Th}\label{lem_Necas}
  There exists a constant $c>0$ such that
  \[
  \|f\|_\Lo\le c(\|f\|_\Hi+\|\nabla f\|_{\Hi^n})
  \]
  for all $f\in\Lo$.
\end{Th}

The following result follows from Theorem \ref{lem_Necas}.

%---------------------------------------------
%---                         Theorem                          ---
%---------------------------------------------
\begin{Th}\label{grad}
  {\rm \cite[Corollary 2.1, 2${}^\circ$]{Girault}}
  There exists a constant $c>0$ such that
  \[
  \|f\|_{\Lo/\R}\le c\|\nabla f\|_{\Hi^n}
  \]
  for all $f\in\Lo$.
\end{Th}

The following two embedding theorems are often called the Poincar\'e inequality.

%---------------------------------------------
%---                         Theorem                          ---
%---------------------------------------------
\begin{Th}\label{P_HR}
  {\rm \cite[Theorem 7.8]{Necas}}
  There exists a constant $c>0$ such that 
  \[
    \|\psi\|_\Lo\le c\|\nabla \psi\|_{\Lo^n}
  \] 
  for all $\psi\in\Hb/\R$.
\end{Th}

%---------------------------------------------
%---                         Theorem                          ---
%---------------------------------------------
\begin{Th}\label{P_HD}
  {\rm \cite[Lemma 3.1]{Girault}}
  There exists a constant $c>0$ such that 
  \[
    \|\psi\|_\Lo\le c\|\nabla \psi\|_{\Lo^n}
  \] 
  for all $\psi\in\HD$.
\end{Th}

%
% ----------------------------------------------
\subsection{Weak formulations of the problems (PP), (S) and (ES)}\label{wform}
% ----------------------------------------------
%

We assume the following conditions for $F,u_b$ and $g_b$:
\begin{align}\label{st_cond}
  F\in\Lo^n,\quad
  u_b\in H^{1/2}(\Gamma),\quad
  \int_\Gamma u_b\cdot \nu=0,
\end{align}
\begin{align}\label{pp_cond}
  g_b\in L^2(\Gamma),\quad
  \dvg F\in\Lo,
\end{align}
\begin{align}\label{Nm_cond}
  \intg g_b=\into\dvg F,
\end{align}
\begin{align}\label{pb_cond}
  p_b\in\Hb.
\end{align}
We start by defining the weak solution to (S).
For all $\varphi\in\Ho^n$, we obtain from the first equation of (S) that
\[
\into F\cdot\varphi
=-\intg \frac{\partial u_\St}{\partial \nu}\cdot\varphi
+\into\nabla u_\St:\nabla\varphi+\into\nabla p_\St\cdot\varphi
=\into\nabla u_\St:\nabla\varphi+\into\nabla p_\St\cdot\varphi.
\]
Using this expression,
the weak form of the Stokes problem becomes as follows:
Find $u_\St\in\Hb^n$ and $p_\St\in\Lo/\R$ such that
\begin{align}\tag{S'}\left\{\begin{array}{ll}
{\DS \into\nabla u_\St:\nabla\varphi+\langle\nabla p_\St,\varphi\rangle
=\into F\cdot\varphi}&\mbox{for all }\varphi\in\Ho^n,\\[8pt]
\dvg u_\St=0&\mbox{in }\Lo,\\
u_\St=u_b&\mbox{in }H^{1/2}(\Gamma)^n.
\end{array}\right.\end{align}

Next, we define the weak formulations of (PP) and (ES)
first for the Neumann boundary condition (\ref{bcn}) 
and them for the mixed boundary condition (\ref{bcm}).
After that, we define generalized weak formulations for (PP) and (ES) 
which cover both cases.

First, we apply the Neumann boundary condition (\ref{bcn}) for (PP) and (ES).
We take a test function $\psi\in\Hb$. From the second equation of (PP), we obtain
\[\begin{array}{ll}
{\DS-\into(\dvg F)\psi}
&{\DS=-\into(\Delta p_\PP)\psi}
={\DS-\intg \frac{\partial p_\PP}{\partial \nu}\psi
+\into\nabla p_\PP\cdot\nabla\psi}\\[8pt]
&{\DS=-\intg g_b\psi+\into\nabla p_\PP\cdot\nabla\psi}.
\end{array}\]
Hence,
\[
\into\nabla p_\PP\cdot\nabla\psi
=\intg g_b\psi-\into(\dvg F)\psi.
\]
We note that $\intg g_b\psi-\into(\dvg F)\psi=\intg g_b[\psi]-\into(\dvg F)[\psi]$
for all $\psi\in\Hb$ by (\ref{Nm_cond}).
Therefore, the weak form of the pressure-Poisson problem
with the Neumann boundary condition (\ref{bcn}) becomes as follows.
Find $u_\PP\in\Hb^n$ and $p_\PP\in\Hb/\R$ such that
\begin{align}\tag{PP$_1$}\left\{\begin{array}{ll}
{\DS \into\nabla u_\PP:\nabla\varphi+\into\nabla p_\PP\cdot\varphi
=\into F\cdot\varphi}&\mbox{for all }\varphi\in\Ho^n,\\[8pt]
{\DS \into\nabla p_\PP\cdot\nabla\psi=\langle G_1,\psi\rangle}&\mbox{for all }\psi\in\Hb/\R,\\
u_\PP=u_b&\mbox{in }H^{1/2}(\Gamma)^n,
\end{array}\right.\end{align}
where $G_1\in\Hb^*$ such that
$\langle G_1,\psi\rangle=\intg g_b\psi-\into(\dvg F)\psi
~(\psi\in\Hb)$.

The weak form of (ES) with the Neumann boundary condition 
can be defined similarly to that of (PP).
Find $u_\eps\in\Hb^n$ and $p_\eps\in\Hb/\R$ such that
\begin{align}\tag{ES$_1$}\left\{\begin{array}{ll}
{\DS \into\nabla u_\eps:\nabla\varphi+\into\nabla p_\eps\cdot\varphi
=\into F\cdot\varphi}&\mbox{for all }\varphi\in\Ho^n,\\[8pt]
{\DS \eps\into\nabla p_\eps\cdot\nabla\psi+\into(\dvg u_\eps)\psi
=\eps\langle G_1,\psi\rangle}&\mbox{for all }\psi\in\Hb/\R,\\
u_\St=u_b&\mbox{in }H^{1/2}(\Gamma)^n.
\end{array}\right.\end{align}

Secondly, we apply the mixed boundary condition (\ref{bcm}) for (PP) and (ES).
We take a test function $\psi\in \HD$. From the second equation of (PP), we obtain
\[\begin{array}{ll}
{\DS-\into(\dvg F)\psi}
&{\DS=-\into(\Delta p_\PP)\psi}
={\DS-\intg \frac{\partial p_\PP}{\partial \nu}\psi
+\into\nabla p_\PP\cdot\nabla\psi}\\[8pt]
&{\DS=-\int_{\Gamma_N} g_b\psi+\into\nabla p_\PP\cdot\nabla\psi}.
\end{array}\]
Hence,
\[
\into\nabla p_\PP\cdot\nabla\psi
=\int_{\Gamma_N} g_b\psi-\into(\dvg F)\psi.
\]
The weak form of the pressure-Poisson problem
with the mixed boundary condition (\ref{bcm}) becomes as follows.
Find $u_\PP\in\Hb^n$ and $p_\PP\in\Hb$ such that
\begin{align}\tag{PP$_2$}\left\{\begin{array}{ll}
{\DS \into\nabla u_\PP:\nabla\varphi+\into\nabla p_\PP\cdot\varphi
=\into F\cdot\varphi}&\mbox{for all }\varphi\in\Ho^n,\\[8pt]
{\DS \into\nabla p_\PP\cdot\nabla\psi=\langle G_2,\psi\rangle}&\mbox{for all }\psi\in\HD,\\
u_\PP=u_b&\mbox{in }H^{1/2}(\Gamma)^n,\\
p_\PP=p_b&\mbox{in }H^{1/2}(\Gamma_D),
\end{array}\right.\end{align}
where $G_2\in\HD^*$ such that
\[
  \langle G_2,\psi\rangle=\int_{\Gamma_N} g_b\psi-\into(\dvg F)\psi
\]
for $\psi\in\HD$.
The weak form of (ES) with the mixed boundary condition (\ref{bcm}) 
can be defined similarly to that of (PP). It reads as follows.
Find $u_\eps\in\Hb^n$ and $p_\eps\in\Hb$ such that
\begin{align}\tag{ES$_2$}\left\{\begin{array}{ll}
{\DS \into\nabla u_\eps:\nabla\varphi+\into\nabla p_\eps\cdot\varphi
=\into F\cdot\varphi}&\mbox{for all }\varphi\in\Ho^n,\\[8pt]
{\DS \eps\into\nabla p_\eps\cdot\nabla\psi+\into(\dvg u_\eps)\psi
=\eps\langle G_2,\psi\rangle}&\mbox{for all }\psi\in\HD,\\
u_\St=u_b&\mbox{in }H^{1/2}(\Gamma)^n,\\
p_\St=p_b&\mbox{in }H^{1/2}(\Gamma_D).
\end{array}\right.\end{align}

Finally, we generalize (PP$_1$) and (PP$_2$) to an abstract pressure-Poisson problem.
Let $Q\subset\Hb$ be a closed subspace as defined in Section \ref{notations}.
Find $u_\PP\in\Hb^n$ and $p_\PP\in Q$ such that
\begin{align}\tag{PP'}
\left\{\begin{array}{ll}
{\DS
\into\nabla u_\PP:\nabla \varphi +\into \nabla p_\PP\cdot\varphi=\into F\cdot\varphi
}& {\rm \mbox{for all }} \varphi\in\Ho^n,\\[8pt]
{\DS
\into\nabla p_\PP\cdot\nabla\psi=\langle G,\psi\rangle
}& {\rm \mbox{for all }} \psi\in Q,\\
u_\PP=u_b &\mbox{in }H^{1/2}(\Gamma)^n,\\
p_\PP-p_b\in Q, &
\end{array}\right.
\end{align}
with $G\in Q^*$.
Indeed, by Theorem \ref{P_HR} and \ref{P_HD}, we obtain (PP$_1$) from (PP') by putting 
$Q:=\Hb/\R$ and $G:=G_1$. Similarly, we obtain (PP$_2$) from (PP') 
by putting $Q:=\HD$ and $G:=G_2$.

We generalize (ES$_1$) and (ES$_2$) to an abstract $\eps$-Stokes problem.
Find $u_\eps\in\Hb^n$ and $p_\eps\in Q$ such that
\begin{align}\tag{ES'}
\left\{\begin{array}{ll}
{\DS
\into\nabla u_\eps:\nabla \varphi +\into \nabla p_\eps\cdot\varphi=\into F\cdot\varphi
}& {\rm \mbox{for all }} \varphi\in\Ho^n,\\[8pt]
{\DS
\eps\into\nabla p_\eps\cdot\nabla\psi+\into(\dvg u_\eps)\psi=\eps\langle G,\psi\rangle
}& {\rm \mbox{for all }} \psi\in Q,\\
u_\eps-u_b\in\Ho^n, &\\
p_\eps-p_b\in Q. &
\end{array}\right.
\end{align}
Indeed, by Theorem \ref{P_HR}, \ref{P_HD}, we obtain (ES$_1$) from (ES') 
by putting $Q:=\Hb/\R$ and $G:=G_1$. Similarly, we also obtain  
(ES$_2$) from (ES') by putting $Q:=\HD$ and $G:=G_2$. 

%
% ----------------------------------------------
\subsection{Well-posedness of (S'), (PP') and (ES')}\label{subsec_well}
% ----------------------------------------------
%

We show the well-posedness of problems (S'), (PP') and (ES')
in Theorem \ref{stokes_t}, \ref{pp_thm} and \ref{estokes_thm}.

%---------------------------------------------
%---                         Theorem                          ---
%---------------------------------------------
\begin{Th}\label{stokes_t}
Under the condition (\ref{st_cond}),
there exists a unique solution 

\noindent$(u_\St,p_\St)\in \Hb^n\times(\Lo/\R)$
satisfying (S').
\end{Th}

%---------------------------------------------
%---                         Proof                              ---
%---------------------------------------------
See \cite[Theorem 2.4 and Remark 2.5]{Temam} for the proof.
%---------------------------------------------

%---------------------------------------------
%---                         Theorem                          ---
%---------------------------------------------
\begin{Th}\label{pp_thm}
Under the condition (\ref{st_cond}) and (\ref{pb_cond}), for $G\in Q^*$,
there exists a unique solution $(u_\PP,p_\PP)\in\Hb^n\times Q$ satisfying (PP').
\end{Th}

%---------------------------------------------
%---                         Proof                              ---
%---------------------------------------------
\noindent $\textit{ Proof.}$
Using the Lax--Milgram theorem,
since $Q\times Q\ni(p,\psi)\mapsto\into\nabla p\cdot\nabla\psi\in\R$
is a continuous and coercive bilinear form,
$p_\PP\in\Hb$ is uniquely determined from the second and fourth equations of (PP').
Then $u_\PP\in\Hb^n$ is also uniquely determined 
from the first and third equations, again using the Lax--Milgram theorem.
\qed%---------------------------------------------

%---------------------------------------------
%---                         Theorem                          ---
%---------------------------------------------
\begin{Th}\label{estokes_thm}
Under the condition (\ref{st_cond}) and (\ref{pb_cond}), 
for $\eps>0$ and $G\in Q^*$,
there exists a unique solution $(u_\eps,p_\eps)\in \Hb^n\times\Hb$
satisfying (ES').
\end{Th}

%---------------------------------------------
%---                         Proof                              ---
%---------------------------------------------
This is a generalization of Theorem 2.6 in \cite{prev}.
See Appendix \ref{prev_thms} for the proof.
%---------------------------------------------

From now on, let the solutions of (S'), (PP') and (ES') be denoted by
$(u_\St,p_\St),(u_\PP,p_\PP)$ and $(u_\eps,p_\eps)$, respectively.
We show their properties in connection with a pressure error 
on the boundary $\Gamma$.

%---------------------------------------------
%---                         Proposition                          ---
%---------------------------------------------
\begin{Prop}\label{sp_prop}
Suppose that $p_\St\in\Hb$, $\Ho\subset Q$ and 
$\langle G,\psi\rangle=$\\ $-\into(\dvg F)\psi$ for all $\psi\in\Ho$.
Then there exists a constant $c>0$ independent of $\eps$ such that
\begin{align}\label{sp_eq}\begin{array}{rrl}
  \|u_\St-u_\PP\|_{\Hb^n}
  &\le& c\|\gamma_0 p_\St-\gamma_0 p_\PP\|_{H^{1/2}(\Gamma)},\\
  \|u_\St-u_\eps\|_{\Hb^n}
  &\le& c\|\gamma_0 p_\St-\gamma_0 p_\PP\|_{H^{1/2}(\Gamma)}.
\end{array}\end{align}
In particular, if $\gamma_0p_\St=\gamma_0 p_\PP$,
then $(u_\St,p_\St)=(u_\PP,p_\PP)=(u_\eps,p_\eps)$ holds for all $\eps>0$.
\end{Prop}

%---------------------------------------------
%---                         Proof                              ---
%---------------------------------------------
This is a generalization of Proposition 2.7 in \cite{prev}.
See Appendix A for the proof.
%---------------------------------------------

Since $\Ho\not\subset\Hb/\R$, Proposition \ref{sp_prop} does not apply directly
for the case of the Neumann boundary condition (\ref{bcn}).
However, we add natural assumptions, then it leads to (\ref{sp_eq}).

%---------------------------------------------
%---                         Proposition                          ---
%---------------------------------------------
\begin{Prop}\label{sp_prop2}
Suppose that $p_\St\in\Hb$, $Q=\Hb/\R$.
If $G\in Q^*=(\Hb/\R)^*$ is such that $G\in\Hb^*$ and
$\langle G,\psi\rangle=\intg g_b\psi-\into(\dvg F)\psi$ 
for all $\psi\in\Hb$, then we have (\ref{sp_eq}).
\end{Prop}

%---------------------------------------------
%---                         Proof                              ---
%---------------------------------------------
\noindent $\textit{ Proof.}$
Since $G\in\Hb^*$ satisfies 
$\langle G,\psi\rangle=\intg g_b\psi-\into(\dvg F)\psi$ for all $\psi\in\Hb$,
it holds that
\[
\into\nabla p_\PP\cdot\nabla\psi=-\into(\dvg F)\psi
\]
for all $\psi\in\Ho$ from the second equation of (PP'). 
Hence, it leads the second equation of (\ref{s-p}).
Using the proof of Proposition \ref{sp_prop}, we obtain (\ref{sp_eq}).
\qed

%
% ----------------------------------------------------------------------------
\section{Links between (ES) and (PP)}\label{sec_ep}
% ----------------------------------------------------------------------------
%

as guaranteed by Theorem \ref{stokes_t}, \ref{pp_thm} 
and \ref{estokes_thm}.
In this section, we show that $(u_\eps,p_\eps)$ converges
to $(u_\PP,p_\PP)$ strongly in $\Hb^n\times\Hb$ 
as $\eps\rightarrow\infty$.
We also treat the case of the regular perturbation asymptotics
by exploring the structure of the lower order terms and their effect
on the convergence rate.

%
% ----------------------------------------------------------------------------
\subsection{Convergence as $\eps\rightarrow\infty$}\label{subsec_ep}
% ----------------------------------------------------------------------------
%

We use the following Lemma \ref{eplem} for the proofs of the theorems 
in this section.

%---------------------------------------------
%---                         Lemma                            ---
%---------------------------------------------
\begin{Lem}\label{eplem}
Let $h\in Q^*$ and $(v_\eps,q_\eps)\in\Ho^n\times Q$ satisfy
  \begin{align}\label{eplem_eq}
  \left\{\begin{array}{ll}
  {\DS
    \into\nabla v_\eps:\nabla \varphi +\into (\nabla q_\eps)\cdot\varphi=0}
  & {\rm \mbox{for all }} \varphi\in\Ho^n,\\[8pt]
  {\DS
    \eps\into\nabla q_\eps\cdot\nabla\psi +\into(\dvg v_\eps)\psi
    =\langle h,\psi\rangle }
  & {\rm \mbox{for all }} \psi \in Q
  \end{array}\right.
  \end{align}
  for an arbitrarily fixed $\eps>0$. Then there exists a constant $c>0$ such that
\[
\|v_\eps\|_{\Hb^n}\le\frac{c}{\eps}\|h\|_{Q^*},\qquad
\|q_\eps\|_\Hb\le\frac{c}{\eps}\|h\|_{Q^*}.
\]
\end{Lem}

%---------------------------------------------
%---                         Proof                              ---
%---------------------------------------------
\noindent $\textit{ Proof.}$
Putting $\varphi:=v_\eps$ and $\psi:=q_\eps$ and
adding two equations of (\ref{eplem_eq}), we obtain
\[
\|\nabla v_\eps\|^2_{\Lo^{n\times n}}
+\eps\|\nabla q_\eps\|^2_{\Lo^n}
\le \|h\|_{Q^*}\|\nabla q_\eps\|_{\Lo^n}.
\]
where we have used
$\into\nabla q_\eps\cdot v_\eps=-\into(\dvg v_\eps)q_\eps$.
Thus
\[
\|\nabla q_\eps\|_{\Lo^n}\le \frac{1}{\eps}\|h\|_{Q^*}.
\]

In addition, from the first equation of (\ref{eplem_eq}) by putting $\varphi:=v_\eps$,
we have
\begin{align*}
  \begin{array}{rl}
  \|\nabla v_\eps\|^2_{\Lo^n}
  = &{\DS \into\nabla v_\eps:\nabla v_\eps}
  =	{\DS-\into(\nabla q_\eps)\cdot v_\eps}
  \le \|\nabla q_\eps\|_{\Lo^n}\|v_\eps\|_{\Lo^n}\\[8pt]
  \le &c\|\nabla q_\eps\|_{\Lo^n}\|\nabla v_\eps\|_{\Lo^{n\times n}},
  \end{array}
\end{align*}and then
\[
\|\nabla v_\eps\|_{\Lo^n}\le c\|\nabla q_\eps\|_{\Lo^n}\le\frac{c}{\eps}\|h\|_{Q^*}.
\]
\qed%---------------------------------------------

Using Lemma \ref{eplem}, we obtain Theorem \ref{ep_conv}.

%---------------------------------------------
%---                         Theorem                          ---
%---------------------------------------------
\begin{Th}\label{ep_conv}
There exists a constant $c>0$ independent of $\eps>0$ such that
\[
\|u_\eps-u_\PP\|_{\Hb^n}\le\frac{c}{\eps}\|\dvg u_\PP\|_{Q^*},\quad
\|p_\eps-p_\PP\|_\Hb\le\frac{c}{\eps}\|\dvg u_\PP\|_{Q^*}.
\]
for all $\eps>0$. In particular, we have
\[
\|u_\eps-u_\PP\|_{\Hb^n}\rightarrow 0,~
\|p_\eps-p_\PP\|_\Hb\rightarrow 0
\quad as~\eps\rightarrow\infty.
\]
\end{Th}

%---------------------------------------------
%---                         Proof                              ---
%---------------------------------------------
\noindent $\textit{ Proof.}$
Combining (PP') and (ES'), we obtain
\begin{align}\label{esp}
  \left\{\begin{array}{ll}
    {\DS
      \into\nabla v_\eps:\nabla\varphi
      +\into \nabla q_\eps\cdot\varphi =0}
    & {\rm \mbox{for all }} \varphi\in\Ho^n,\\[8pt]
    {\DS
      \eps\into\nabla q_\eps\cdot\nabla\psi+\into(\dvg v_\eps)\psi
	=-\into(\dvg u_\PP)\psi}
    & {\rm \mbox{for all }} \psi \in Q,
  \end{array}\right.
\end{align}
where $v_\eps:=u_\eps-u_\PP,q_\eps:=p_\eps-p_\PP$ and $h:=\dvg u_\PP$.
By Lemma \ref{eplem}, we conclude the proof.
\qed%---------------------------------------------

%---------------------------------------------
%---                         Corollary                          ---
%---------------------------------------------
\begin{Cor}\label{dvgup}
If $u_\PP$ satisfies $\dvg u_\PP=0$,
then $u_\eps=u_\PP$ and $p_\eps=p_\PP$
hold for all $\eps>0$. Furthermore,
$u_\St=u_\eps=u_\PP$ and $p_\St=[p_\eps]=[p_\PP]$ hold for all $\eps>0$.
\end{Cor}

%
% ----------------------------------------------
\subsection{Regular Perturbation Asymptotics}\label{sec_RPA}
% ----------------------------------------------
%

%\[
%u_\eps=u_\PP+\frac{1}{\eps}v_{(1)}+\cdots
%\]
%Rem. Idea of operator (Neumann)

By Theorem \ref{ep_conv}, we have that 
$\|\eps(u_\eps-u_\PP)\|_{\Hb^n}\le c$ and $\|\eps(p_\eps-p_\PP)\|_\Hb\le c$
for all $\eps>0$.
It implies that there exists a subsequence of
$(\eps(u_\eps-u_\PP),\eps(p_\eps-p_\PP))$ which converges weakly to
$(v^{(1)},q^{(1)})\in\Ho^n\times Q$ if $\eps\rightarrow\infty$.
The next theorem states properties of the limit functions $v^{(1)}$ and $q^{(1)}$.

%---------------------------------------------
%---                         Theorem                          ---
%---------------------------------------------
\begin{Th}\label{v1t}
  Let $v^{(1)}_\eps:=\eps(u_\eps-u_\PP)\in\Ho^n,q^{(1)}_\eps:=\eps(p_\eps-p_\PP)\in Q$
  and let $(v^{(1)},q^{(1)})\in\Ho^n\times Q$ satisfy
\begin{align}\label{v1eq}
  \left\{
  \begin{array}{ll}{\displaystyle
  \into \nabla v^{(1)}:\nabla \varphi+\into (\nabla q^{(1)})\cdot\varphi=0}
    & {\rm \mbox{for all }} \varphi \in\Ho^n, \\[8pt]
    {\displaystyle
    \into \nabla q^{(1)}\cdot \nabla \psi = -\into(\dvg u_\PP)\psi }
    & {\rm \mbox{for all }} \psi\in Q.
  \end{array}
  \right.
\end{align}
Then there exists a constant $c>0$ independent of $\eps$ such that
\begin{align*}
  \|v^{(1)}_\eps-v^{(1)}\|_{\Hb^n}\le\frac{c}{\eps}\|\dvg v^{(1)}\|_{Q^*},\qquad
  \|q^{(1)}_\eps-q^{(1)}\|_\Hb\le\frac{c}{\eps}\|\dvg v^{(1)}\|_{Q^*}.
\end{align*}
\end{Th}

%---------------------------------------------
%---                         Proof                              ---
%---------------------------------------------
\noindent $\textit{ Proof.}$
The existence and the uniqueness of the pair $(v^{(1)},q^{(1)})\in\Ho^n\times Q$
as a solution to (\ref{v1eq}) follows from Theorem \ref{pp_thm}. 
As in (\ref{esp}), we have
\begin{align}\label{v1e}
  \left\{
  \begin{array}{ll}
    {\displaystyle
      \into\nabla v^{(1)}_\eps:\nabla\varphi +\into (\nabla q^{(1)}_\eps)\cdot\varphi =0  }
    & {\rm \mbox{for all }} \varphi\in\Ho^n,\\[8pt]
    {\displaystyle
      \into\nabla q^{(1)}_\eps\cdot\nabla\psi+\frac{1}{\eps}\into(\dvg v^{(1)}_\eps)\psi=
      -\into (\dvg u_\PP)\psi}
    & {\rm \mbox{for all }} \psi \in Q.
  \end{array}
  \right.
\end{align}
Subtracting (\ref{v1eq}) from (\ref{v1e}), it holds that
\begin{align*}
  \left\{
  \begin{array}{ll}
    {\displaystyle
      \into\nabla (v^{(1)}_\eps-v^{(1)}):\nabla\varphi
      +\into (\nabla (q^{(1)}_\eps-q^{(1)}))\cdot\varphi =0  }
    & {\rm \mbox{for all }} \varphi\in\Ho^n,\\[8pt]
    {\displaystyle
      \into\nabla (q^{(1)}_\eps-q^{(1)})\cdot\nabla\psi
      +\frac{1}{\eps}\into(\dvg v^{(1)}_\eps)\psi=0}
    & {\rm \mbox{for all }} \psi \in Q.
  \end{array}
  \right.
\end{align*}
Hence,
\[
  \left\{\begin{array}{ll}
    {\displaystyle
      \into\nabla v_\eps:\nabla\varphi
      +\into \nabla q_\eps\cdot\varphi =0  }
    & {\rm \mbox{for all }} \varphi\in\Ho^n,\\[8pt]
    {\displaystyle
      \eps\into\nabla q_\eps\cdot\nabla\psi
      +\into(\dvg v_\eps)\psi=-\into(\dvg v^{(1)})\psi}
    & {\rm \mbox{for all }} \psi \in Q,
  \end{array}\right.
\]
where $v_\eps:=v^{(1)}_\eps-v^{(1)}, q_\eps:=q^{(1)}_\eps-q^{(1)}$ and $h:=-\dvg v^{(1)}$.
By Lemma \ref{eplem} , we have
$$
\|v^{(1)}_\eps-v^{(1)}\|_{\Hb^n}\le\frac{c}{\eps}\|\dvg v^{(1)}\|_{Q^*},\quad
\|q^{(1)}_\eps-q^{(1)}\|_\Hb\le \frac{c}{\eps}\|\dvg v^{(1)}\|_{Q^*}
$$
for all $\eps>0$.
\qed%---------------------------------------------

Next, we generalize Theorem \ref{v1t} to the following theorem:

%---------------------------------------------
%---                         Theorem                          ---
%---------------------------------------------
\begin{Th}\label{reg}
  Let $k\in {\mathbb N}$ be arbitrary ($k\ge 1$) and let $v^{(0)}:=u_\PP$.
  If functions $v^{(1)},v^{(2)},\cdots,v^{(k)}\in\Ho^n$ and
  $q^{(1)},q^{(2)},\cdots,q^{(k)}\in Q$ satisfy
  \begin{align}\label{vqi}
  \left\{
  \begin{array}{ll}{\displaystyle
  \into \nabla v^{(i)}:\nabla \varphi+\into (\nabla q^{(i)})\cdot\varphi=0}
    & {\rm \mbox{for all }} \varphi \in\Ho^n,\\[8pt]
    {\displaystyle
    \into \nabla q^{(i)}\cdot \nabla \psi = -\into(\dvg  v^{(i-1)})\psi }
    & {\rm \mbox{for all }} \psi\in Q,
  \end{array}
  \right.
  \end{align}
  for all $1\le i\le k$, then there exists a constant $c>0$
independent of $\eps$ satisfying
  \begin{equation*}
    \begin{array}{c}
      {\displaystyle \left\|u_\eps
          -\left(u_\PP+\frac{1}{\eps}v^{(1)}
          +\cdots+\left(\frac{1}{\eps}\right)^{k}v^{(k)}\right)\right\|_{\Hb^n}
	\le\frac{c}{\eps^{k+1}}\|\dvg  v^{(k)}\|_{Q^*},}\\[8pt]
      {\displaystyle \left\|p_\eps
          -\left(p_\PP+\frac{1}{\eps}q^{(1)}
          +\cdots+\left(\frac{1}{\eps}\right)^{k}q^{(k)}\right)\right\|_\Hb
	\le\frac{c}{\eps^{k+1}}\|\dvg  v^{(k)}\|_{Q^*}.}
    \end{array}
  \end{equation*}
\end{Th}

%---------------------------------------------
%---                         Proof                              ---
%---------------------------------------------
\noindent $\textit{ Proof.}$
Let $(v^{(i)}_\eps,q^{(i)}_\eps)\in\Ho^n\times Q~(1\le i\le k)$ satisfy
\begin{align}\label{vqei}
  \left\{
  \begin{array}{ll}{\displaystyle
      \into\nabla v^{(i)}_\eps:\nabla\varphi+\into(\nabla q^{(i)}_\eps)\cdot\varphi=0}
    & {\rm \mbox{for all }} \varphi \in\Ho^n,\\[8pt]
    {\displaystyle
      \into\nabla q^{(i)}_\eps\cdot\nabla\psi
      +\frac{1}{\eps}\into(\dvg v^{(i)}_\eps)\psi
      = -\into(\dvg  v^{(i-1)})\psi }
    & {\rm \mbox{for all }} \psi\in Q.
  \end{array}
  \right.
\end{align}
Subtracting (\ref{vqi}) from (\ref{vqei}), it holds that
\begin{align*}
  \left\{
  \begin{array}{ll}
    {\displaystyle
      \into\nabla (v^{(i)}_\eps-v^{(i)}):\nabla\varphi
      +\into (\nabla (q^{(i)}_\eps-q^{(i)}))\cdot\varphi =0  }
    & {\rm \mbox{for all }} \varphi\in\Ho^n,\\[8pt]
    {\displaystyle
      \into\nabla (q^{(i)}_\eps-q^{(i)})\cdot\nabla\psi
      +\frac{1}{\eps}\into(\dvg v^{(i)}_\eps)\psi=0}
    & {\rm \mbox{for all }} \psi \in Q.
  \end{array}
  \right.
\end{align*}
Setting $v_\eps:=v^{(i)}_\eps-v^{(i)},q_\eps:=q^{(i)}_\eps-q^{(i)}$ and 
$h:=-\dvg v^{(i)}$, we obtain from Lemma \ref{eplem} that the estimates
\[
\|v^{(i)}_\eps-v^{(i)}\|_{\Hb^n}\le\frac{c}{\eps}\|\dvg  v^{(i)}\|_{Q^*},\quad
\|q^{(i)}_\eps-q^{(i)}\|_\Hb\le \frac{c}{\eps}\|\dvg  v^{(i)}\|_{Q^*}
\]
hold for all $\eps>0$.
In particular, putting $i:=k$, we obtain
\[
\|v^{(k)}_\eps-v^{(k)}\|_{\Hb^n}\le\frac{c}{\eps}\|\dvg  v^{(k)}\|_{Q^*},
\]
\[
\|q^{(k)}_\eps-q^{(k)}\|_\Hb\le \frac{c}{\eps}\|\dvg  v^{(k)}\|_{Q^*}
\]
for all $\eps>0$.
By the uniqueness of the solution to (ES') in Theorem \ref{estokes_thm},
it leads that 
$v^{(i+1)}_\eps=\eps(v^{(i)}_\eps-v^{(i)}),
q^{(i+1)}_\eps=\eps(q^{(i)}_\eps-q^{(i)})$ for all $i=1,\cdots,k-1$,
and thus
\begin{align*}
  \begin{array}{rl}
  &v^{(k)}_\eps-v^{(k)}\\[8pt]
  =&\eps(v^{(k-1)}_\eps-v^{(k-1)})-v^{(k)}\\[8pt]
  =&{\DS\eps \left(v^{(k-1)}_\eps-\left(
  v^{(k-1)}+\left(\frac{1}{\eps}\right)v^{(k)}\right)\right)}\\[8pt]
  =&\cdots\\[8pt]
  =&{\DS\eps^{k-1}\left(v^{(1)}_\eps-\left(
  v^{(1)}+\cdots+ \left(\frac{1}{\eps}\right)^{k-2}v^{(k-1)}
  +\left(\frac{1}{\eps}\right)^{k-1}v^{(k)}\right)\right)}\\[8pt]
  =&{\DS\eps^{k}\left(u_\eps-\left(
  u_\PP+\frac{1}{\eps}v^{(1)}
  +\cdots+ \left(\frac{1}{\eps}\right)^{k-1}v^{(k-1)}
  +\left(\frac{1}{\eps}\right)^{k}v^{(k)}\right)\right)},
  \end{array}
\end{align*}
\begin{align*}
  \begin{array}{rl}
  &q^{(k)}_\eps-q^{(k)}\\[8pt]
  =&\eps(q^{(k-1)}_\eps-q^{(k-1)})-q^{(k)}\\[8pt]
  =&{\DS\eps \left(q^{(k-1)}_\eps-\left(
  q^{(k-1)}+\left(\frac{1}{\eps}\right)q^{(k)}\right)\right)}\\[8pt]
  =&\cdots\\[8pt]
  =&{\DS\eps^{k-1}\left(q^{(1)}_\eps-\left(
  q^{(1)}+\cdots+ \left(\frac{1}{\eps}\right)^{k-2}q^{(k-1)}
  +\left(\frac{1}{\eps}\right)^{k-1}q^{(k)}\right)\right)}\\[8pt]
  =&{\DS\eps^{k}\left(p_\eps-\left(
  p_\PP+\frac{1}{\eps}q^{(1)}
  +\cdots+ \left(\frac{1}{\eps}\right)^{k-1}q^{(k-1)}
  +\left(\frac{1}{\eps}\right)^{k}q^{(k)}\right)\right)}.
  \end{array}
\end{align*}
Hence it holds that
\begin{align*}
  \left\|u_\eps-\left(
  u_\PP+\frac{1}{\eps}v^{(1)}
  +\cdots+\left(\frac{1}{\eps}\right)^{k}v^{(k)}\right)\right\|_{\Hb^n}
  &\le\frac{c}{\eps^{k+1}}\|\dvg  v^{(k)}\|_{Q^*},\\[8pt]
  \left\|p_\eps-\left(
  p_\PP+\frac{1}{\eps}q^{(1)}
  +\cdots+\left(\frac{1}{\eps}\right)^{k}q^{(k)}\right)\right\|_\Hb
  &\le\frac{c}{\eps^{k+1}}\|\dvg  v^{(k)}\|_{Q^*}.
\end{align*}
\qed%---------------------------------------------

%---------------------------------------------
%---                         Remark                          ---
%---------------------------------------------
\begin{Rem}
Theorem \ref{reg} can be interpreted from the operator theory.

Let $t\ge0,X:=\Ho^n\times Q,Y:=\Hi^n\times Q^*$ be equipped with norms
\[\begin{array}{rl}
  \|(u,p)\|^2_X&:=\|u\|^2_{\Hb^n}+\|p\|^2_\Hb,\\[8pt]
  \|(f,g)\|^2_Y&:=\|f\|^2_{\Hi^n}+\|g\|^2_{Q^*}
\end{array}\]
for $(u,p)\in X,(f,g)\in Y$, and let $A$ and $B$ be
\[\begin{array}{rccc}
A:&X&\longrightarrow&Y\\
&\rotatebox{90}{$\in$}&&\rotatebox{90}{$\in$}\\
&(u,p)&\longmapsto&
(-\Delta u+\nabla p, \Delta p),
\end{array}\quad\begin{array}{rccc}
B:&X&\longrightarrow&Y\\
&\rotatebox{90}{$\in$}&&\rotatebox{90}{$\in$}\\
&(u,p)&\longmapsto&
(0, \dvg u).
\end{array}\]
Then $(u_\PP,p_\PP)$ and $(u_\eps,p_\eps)$ satisfy
\[
  A(u_\PP,p_\PP)=f,\quad \left(A+\frac{1}{\eps}B\right)(u_\eps,p_\eps)=f,
\]
where $f=(F,G)$.
We have $A+tB\in{\rm Isom}(X,Y)$ for an arbitrary $t\ge 0$
by the analogy of Theorem \ref{pp_thm} ($t=0$)
and Theorem \ref{estokes_thm} ($t=1/\eps$).
Equation (\ref{vqi}) states that 
\[
  A(v^{(i)},q^{(i)})=-B(v^{(i-1)},q^{(i-1)})
\]
for $i=1,\cdots,k$, i.e.
\[\begin{array}{rl}
  (v^{(i)},q^{(i)})&=-A^{-1} B(v^{(i-1)},q^{(i-1)})
  =\cdots=(-A^{-1} B)^i(u_\PP,p_\PP)\\[8pt]
  &=A^{-1}(-BA^{-1})^i f.
\end{array}\]
By Theorem \ref{reg}, there exists a constant $c>0$ such that
\[
\left\|\left(A+\frac{1}{\eps}B\right)^{-1} f
-A^{-1}\sum^{k}_{i=0}\left(-\frac{1}{\eps}BA^{-1}\right)^i f\right\|_X
\le \frac{c}{\eps^{k+1}}\|(BA^{-1})^{k+1} f\|_Y
\]
for all $\eps>0,f\in Y$.
\end{Rem}

%
% ----------------------------------------------------------------------------
\section{Convergence of (ES) to (S)}\label{sec_es}
% ----------------------------------------------------------------------------
%

In this section, we show that $(u_\eps,p_\eps)$ converges
to $(u_\St,p_\St)$ weakly in $\Ho^n\times(\Lo/\R)$ 
as $\eps\rightarrow 0$.
Moreover, if $p_\St\in\Hb$, then $(u_\eps,p_\eps)$ converges
to $(u_\St,p_\St)$ strongly in $\Ho^n\times(\Lo/\R)$ 
as $\eps\rightarrow 0$.

The outline of the proof of our convergence results 
(Theorem \ref{bb}, \ref{es_conv} and \ref{esn_conv}) is as follows.
First, we prove the boundedness of the sequence $((u_\eps,p_\eps))_{\eps>0}$
in $\Ho^n\times(\Lo/\R)$.
By the reflexivity of $\Ho^n\times(\Lo/\R)$, the sequence has a subsequence
converging weakly in $\Ho^n\times(\Lo/\R)$.
In the end, we show that the limit pair of functions satisfies (S').

We start this section with a useful lemma. 

%---------------------------------------------
%---                         Lemma                            ---
%---------------------------------------------
\begin{Lem}\label{pleu}
  If $v\in\Hb^n,q\in\Lo$ and $f\in\Hi^n$ satisfy
  \[
  \into\nabla v:\nabla\varphi +\langle\nabla q,\varphi\rangle=\langle f,\varphi\rangle
  \quad \mbox{for all }\varphi\in\Ho^n,
  \]
  then there exists a constant $c>0$ such that
  \[
  \|q\|_{\Lo/\R}\le c(\|\nabla v\|_{\Lo^{n\times n}}+\|f\|_{\Hi^n}).
  \]
\end{Lem}

%---------------------------------------------
%---                         Proof                              ---
%---------------------------------------------
\noindent $\textit{ Proof.}$
Let $c$ be the constant from Theorem \ref{grad}. Then we obtain
\begin{align*}
  \begin{array}{rl}
  \|q\|_{\Lo/\R}
  \le&c\|\nabla q\|_{\Hi^n}
  ={\DS c\sup_{\varphi\in S_n}
  |\langle\nabla q,\varphi\rangle|}\\[8pt]
  \le&{\DS c\sup_{\varphi\in S_n}
  \left(\left|\into \nabla v:\nabla\varphi\right|
    +|\langle f,\varphi\rangle|\right)}\\[8pt]
  \le&c(\|\nabla v\|_{\Lo^{n\times n}}+\|f\|_{\Hi^n}).
  \end{array}
\end{align*}
\qed%---------------------------------------------

%---------------------------------------------
%---                         Theorem                          ---
%---------------------------------------------
\begin{Th}\label{bb}
There exists a constant $c>0$ independent of $\eps$ such that
\[
\|u_\eps\|_{\Hb^n}\le c,\quad\|p_\eps\|_{\Lo/\R}\le c
\quad{\mbox{for all }\eps>0}.
\]
Furthermore, if $C^\infty_0(\Omega)\subset Q$,
then we obtain
\[
u_{\eps}\rightharpoonup u_\St~weakly~\mbox{in }\Hb^n,~
[p_{\eps}]\rightharpoonup p_\St~weakly~\mbox{in }\Lo/\R
\quad{\rm as}~\eps\rightarrow 0.
\]
\end{Th}

%---------------------------------------------
%---                         Proof                              ---
%---------------------------------------------
See Appendix A for the proof.

If we add a regularity assumption of $p_\St$, 
then $(u_\eps,p_\eps)$ converges strongly in $\Hb^n\times\Lo/\R$
%---------------------------------------------
%---                         Theorem                          ---
%---------------------------------------------
\begin{Th}\label{es_conv}
Suppose that $p_\St\in\Hb$. Then we obtain
\[
u_{\eps}\rightarrow u_\St~strongly~\mbox{in }\Hb^n,~
[p_{\eps}]\rightarrow p_\St~strongly~\mbox{in }\Lo/\R
\quad{\rm as}~\eps\rightarrow 0.
\]
\end{Th}

%---------------------------------------------
%---                         Proof                              ---
%---------------------------------------------
See Appendix A for the proof.

Theorem \ref{es_conv} does not give the convergence rate.
If $Q=\Hb/\R$ (corresponding to the Neumann boundary condition (\ref{bcn})), 
then the convergence rate becomes $\sqrt{\eps}$.

%---------------------------------------------
%---                         Theorem                          ---
%---------------------------------------------
\begin{Th}\label{esn_conv}
Suppose that $Q=\Hb/\R$ and $p_\St\in\Hb$. Then there exists a constant $c>0$
independent of $\eps$ such that
\[
\|u_\eps-u_\St\|_{\Hb^n}\le c\sqrt\eps,~
\|p_\eps-p_\St\|_{\Lo}\le c\sqrt\eps.
\]
\end{Th}

%---------------------------------------------
%---                         Proof                              ---
%---------------------------------------------
\noindent $\textit{ Proof.}$
We obtain from (ES') and (S') that
\[
  \left\{
  \begin{array}{ll}
    {\displaystyle
  \into\nabla (u_\eps-u_\St):\nabla\varphi +\into (\nabla(p_\eps-p_\St))\cdot\varphi =0  }
    & {\rm \mbox{for all }} \varphi\in\Ho^n,\\[8pt]
    {\displaystyle
    \eps\into\nabla p_\eps\cdot\nabla\psi+\into(\dvg u_\eps)\psi=\eps\langle G,\psi\rangle}
    & {\rm \mbox{for all }} \psi \in\Hb/\R.
  \end{array}
  \right.
\]
Putting $\varphi:=u_\eps-u_\St\in\Ho^n$ and $\psi:=p_\eps-p_\St\in\Hb/\R$, we get
\begin{align}\label{es_mid}\begin{array}{ll}
&\|\nabla(u_\eps-u_\St)\|^2_{\Lo^{n\times n}}
+\eps\into\nabla p_\eps\cdot\nabla(p_\eps-p_\St)\\[8pt]
=&{\DS -\into (\nabla(p_\eps-p_\St))\cdot(u_\eps-u_\St)
-\into(\dvg u_\eps)(p_\eps-p_\St)
+\eps\langle G,p_\eps-p_\St\rangle}\\[8pt]
=&{\DS \into (\dvg u_\eps-\dvg u_\St)(p_\eps-p_\St)
-\into(\dvg u_\eps)(p_\eps-p_\St)
+\eps\langle G,p_\eps-p_\St\rangle}\\[8pt]
=&{\DS \eps\langle G,p_\eps-p_\St\rangle.}\\[8pt]
\end{array}\end{align}
Subtracting $\eps\into\nabla p_\St\cdot\nabla(p_\eps-p_\St)$ 
from both sides of (\ref{es_mid}), we obtain
\begin{align}\label{eq_esn}\begin{array}{ll}
&\|\nabla(u_\eps-u_\St)\|^2_{\Lo^{n\times n}}+\eps\|\nabla (p_\eps-p_\St)\|^2_{\Lo^n}\\[8pt]
=&{\DS -\eps\into\nabla p_\St\cdot\nabla(p_\eps-p_\St)
+\eps\langle G,p_\eps-p_\St\rangle}\\[8pt]
\le &\eps(\|\nabla p_\St\|_{\Lo^n}+\|G\|_{(\Hb/\R)^*})\|\nabla(p_\eps-p_\St)\|_{\Lo^n}.
\end{array}\end{align}
To clarify the following estimates, we set 
$\alpha:=\|\nabla(u_\eps-u_\St)\|_{\Lo^{n\times n}},
\beta:=\|\nabla(p_\eps-p_\St)\|_{\Lo^n},a:=\|\nabla p_\St\|_{\Lo^n}+\|G\|_{(\Hb/\R)^*}$.
The estimate (\ref{eq_esn}) reads as
\[
\alpha^2+\eps \beta^2\le\eps a\beta,~
\left(\frac{\alpha}{\sqrt{\eps}}\right)^2+\left(\beta-\frac{a}{2}\right)^2
\le\left(\frac{a}{2}\right)^2.
\]
Hence, $\alpha\le a\sqrt{\eps}/2$, i.e.,
$\|\nabla(u_\eps-u_\St)\|_{\Lo^{n\times n}}\le
 (\sqrt{\eps}/2)(\|\nabla p_\St\|_{\Lo^n}+\|G\|_{(\Hb/\R)^*})$.
By Lemma \ref{pleu}, we obtain
\[\begin{array}{ll}
\|p_\eps-p_\St\|_\Lo&\le c\|\nabla(u_{\eps_k}-u_\St)\|_{\Lo^{n\times n}}
=c\alpha
\le{\DS c\frac{a\sqrt{\eps}}{2}}\\[8pt]
&={\DS c\frac{\sqrt{\eps}}{2}(\|\nabla p_\St\|_{\Lo^n}+\|G\|_{(\Hb/\R)^*}).}
\end{array}\]
\qed%---------------------------------------------

%
% ----------------------------------------------------------------------------
\section{Numerical examples}\label{sec_numex}
% ----------------------------------------------------------------------------
%

%stabilized method, example of numerical simulation.

\begin{figure}[H]
  \begin{minipage}{0.47\hsize}\centering
  \includegraphics[width=\linewidth,bb=50 50 410 302]{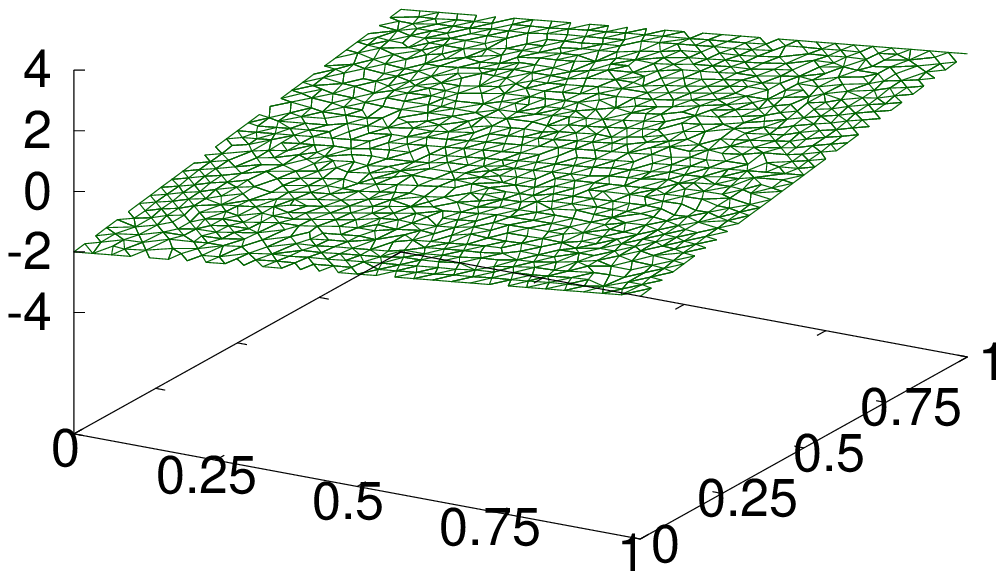}
    \end{minipage}
  \centering\begin{minipage}{0.47\hsize}\centering
  \includegraphics[width=\linewidth,bb=50 50 410 302]{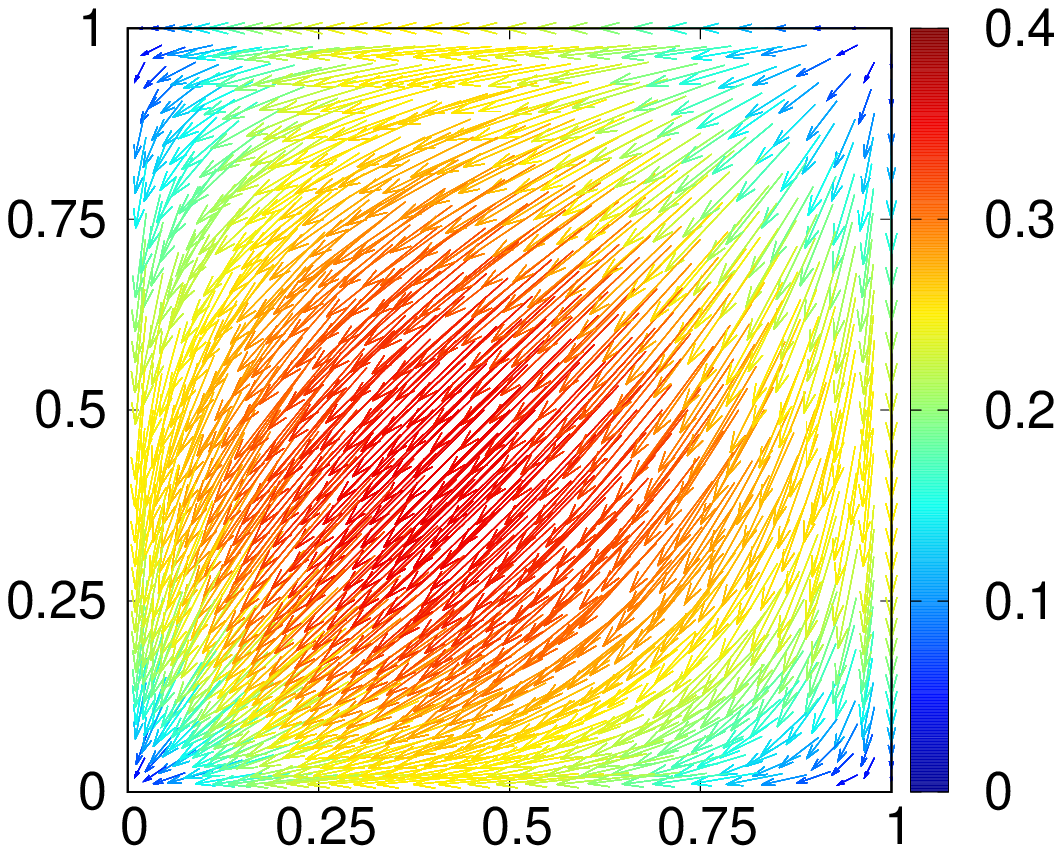}
  \end{minipage}
  \caption{$p_\PP$ (left) and $u_\PP$ (right).
  The color scale indicates the length of $|u_\PP(\xi)|$
  at each node $\xi$.}\label{pp_fig}
\end{figure}

\begin{figure}[H]
\centering\begin{minipage}{0.47\hsize}\centering
\includegraphics[width=\linewidth,bb=50 50 410 302]{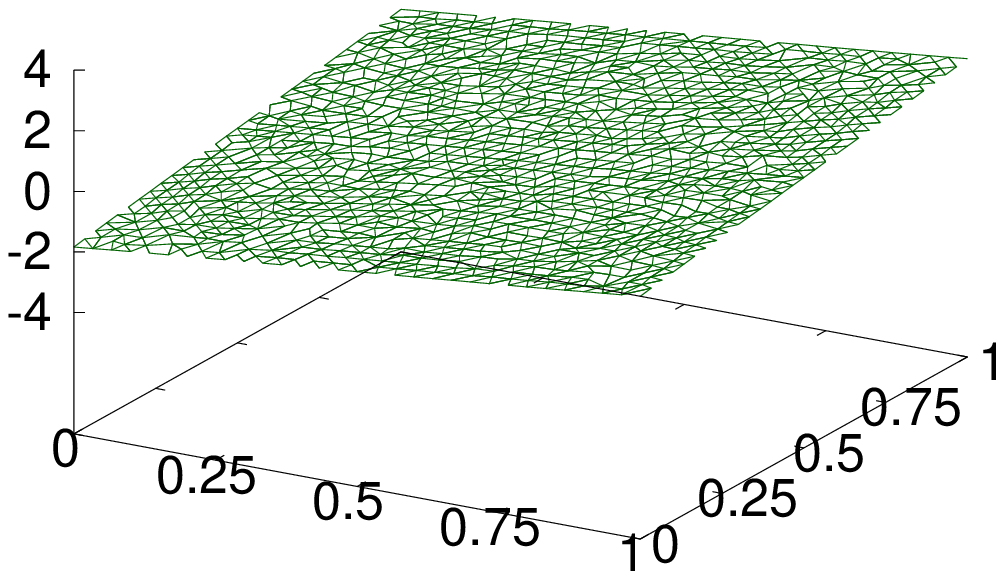}
\subcaption{}\label{p0_fig}
\end{minipage}
\begin{minipage}{0.47\hsize}\centering
\includegraphics[width=\linewidth,bb=50 50 410 302]{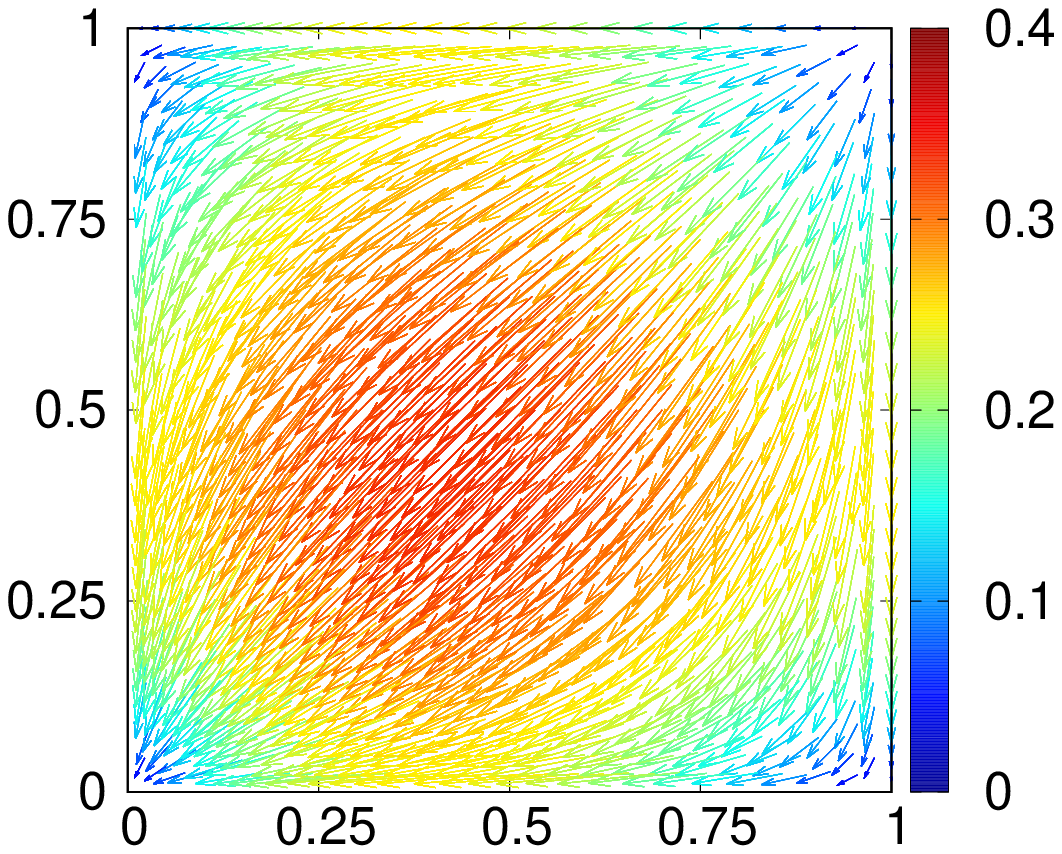}
\subcaption{}\label{u0_fig}
\end{minipage}
\centering\begin{minipage}{0.47\hsize}\centering
\includegraphics[width=\linewidth,bb=50 50 410 302]{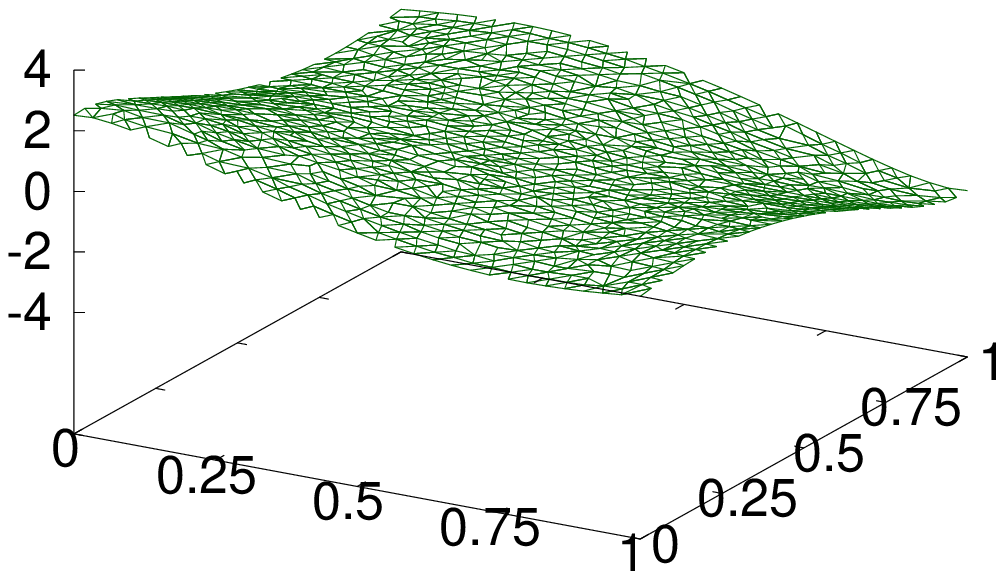}
\subcaption{}\label{p2_fig}
\end{minipage}
\begin{minipage}{0.47\hsize}\centering
\includegraphics[width=\linewidth,bb=50 50 410 302]{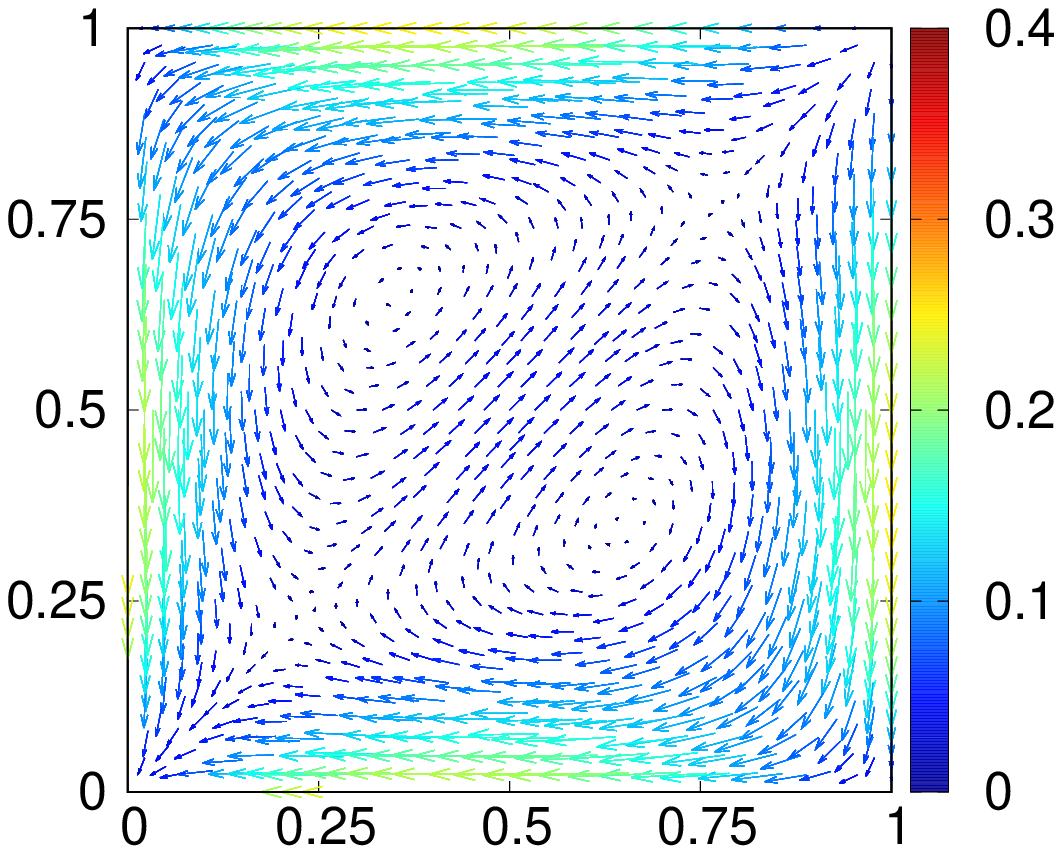}
\subcaption{}\label{u2_fig}
\end{minipage}
\centering\begin{minipage}{0.47\hsize}\centering
\includegraphics[width=\linewidth,bb=50 50 410 302]{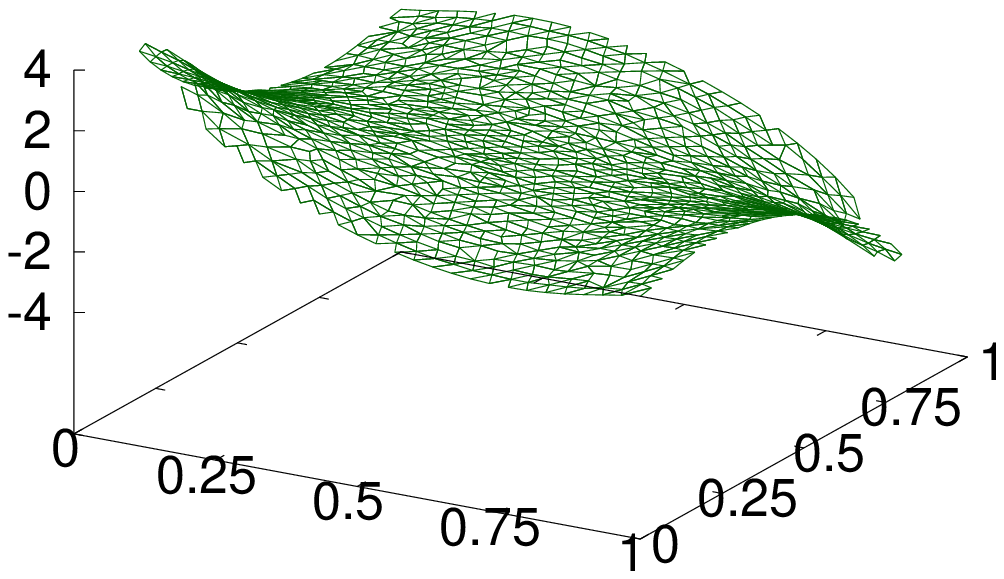}
\subcaption{}\label{p4_fig}
\end{minipage}
\begin{minipage}{0.47\hsize}\centering
\includegraphics[width=\linewidth,bb=50 50 410 302]{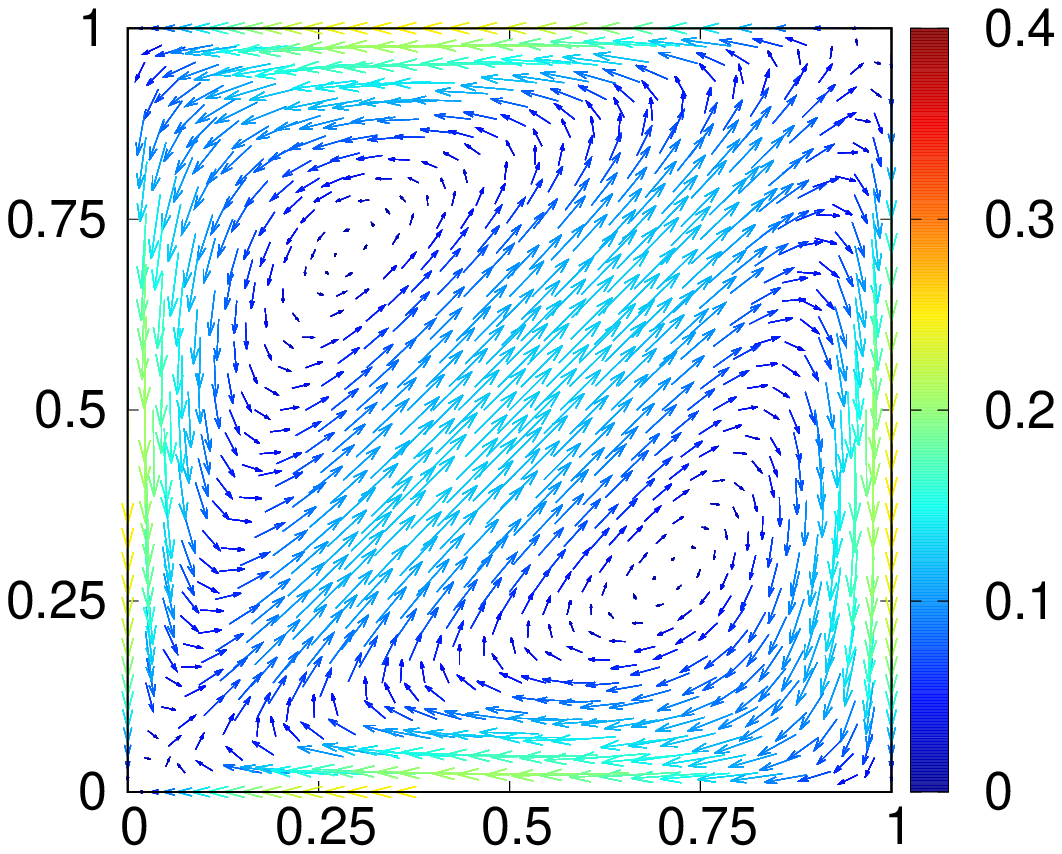}
\subcaption{}\label{u4_fig}
\end{minipage}
\caption{
$p_\eps$ (a) and $u_\eps$ (b) with $\eps=1$.
$p_\eps$ (c) and $u_\eps$ (d) with $\eps=10^{-2}$.
$p_\eps$ (e) and $u_\eps$ (f) with $\eps=10^{-4}$.
The color scales indicate the length of $|u_\eps(\xi)|$
at each node $\xi$.}\label{es_fig}
\end{figure}

\begin{figure}[H]
  \begin{minipage}{0.47\hsize}\centering
  \includegraphics[width=\linewidth,bb=50 50 410 302]{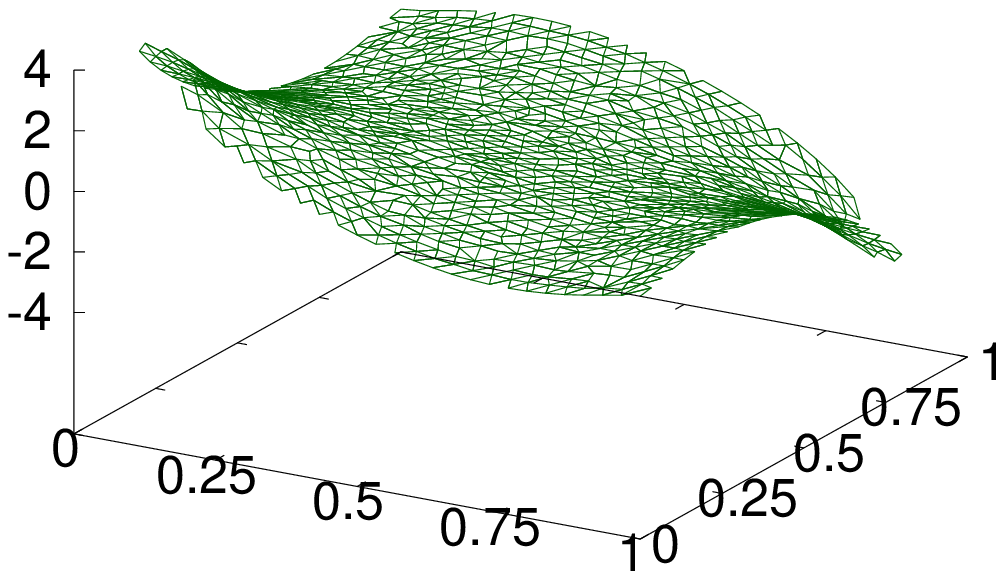}
    \end{minipage}
  \centering\begin{minipage}{0.47\hsize}\centering
  \includegraphics[width=\linewidth,bb=50 50 410 302]{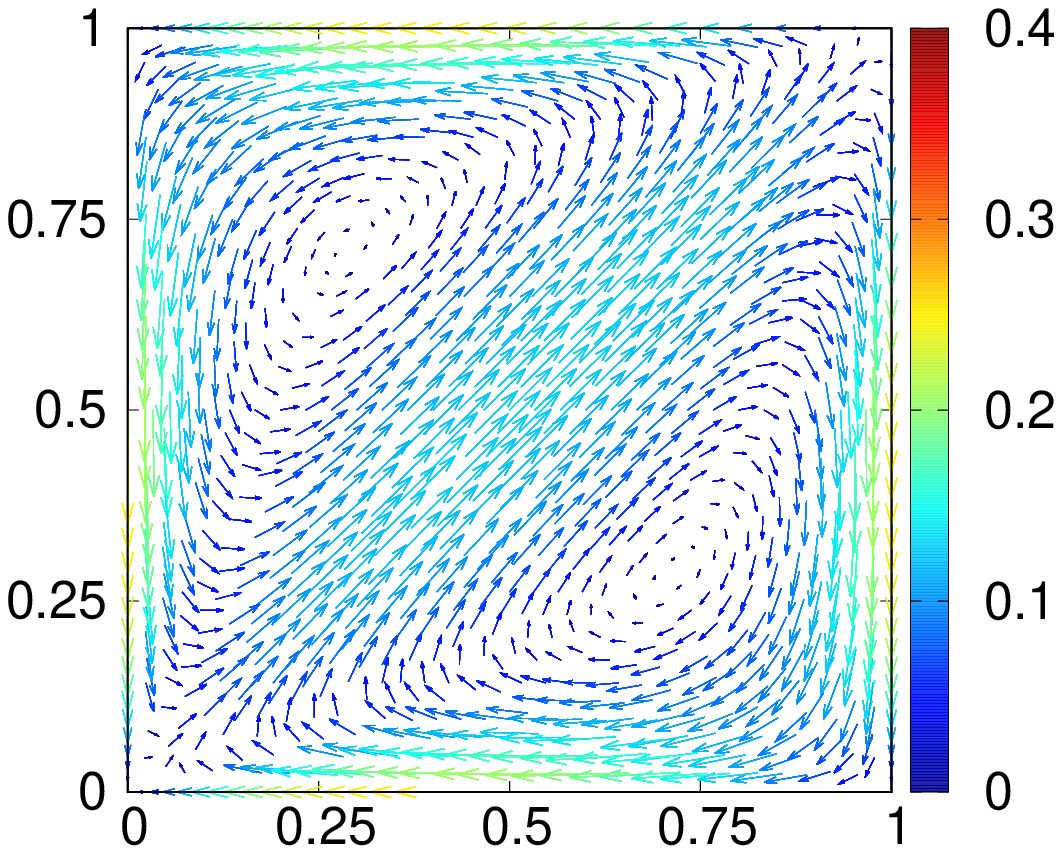}
  \end{minipage}
  \caption{$p_\St$ (left) and $u_\St$ (right).
  The color scale indicates the length of $|u_\St(\xi)|$
  at each node $\xi$.}\label{s_fig}
\end{figure}

%If $\eps$ is small enough, (ES) can be solved by the stabilized finite
%element method with P1/P1 element as a substitute for (S).
For our simulations, we consider $\Omega=(0,1)\times (0,1)$.
We take the following boundary conditions:
\[
u_b=(x(x-1),y(y-1))^T,~g_b=(2,2)^T\cdot\nu
\]
on $\Gamma$. The exact solutions for (PP$_1$) are 
$u_\PP=(x(x-1),y(y-1))^T$ and $p_\PP=2x+2y-2$.
We solve the problems (PP$_1$), (ES$_1$) and (S') numerically by
using the finite element method with P2/P1 elements 
by the software FreeFem++ \cite{FreeFem}. The numerical solutions
$(u_\PP,p_\PP),(u_\eps,p_\eps)~(\eps=1,10^{-2}{\rm~or~}10^{-4})$
and $(u_\St,p_\St)$ to the problems (PP$_1$), (ES$_1$) and (S'), respectively,
are illustrated in Fig. \ref{pp_fig}--\ref{s_fig}.
From these pictures we observe that $(u_\eps,p_\eps)$ seems to 
converge to $(u_\PP,p_\PP)$ as $\eps\rightarrow\infty$ and
to $(u_\St,p_\St)$ as $\eps\rightarrow 0$ 
(as expected from Theorem \ref{ep_conv} and \ref{es_conv}.) 

\begin{figure}[H]
\centering\begin{minipage}{0.47\hsize}\centering
\includegraphics[width=\linewidth,bb=50 50 410 302]{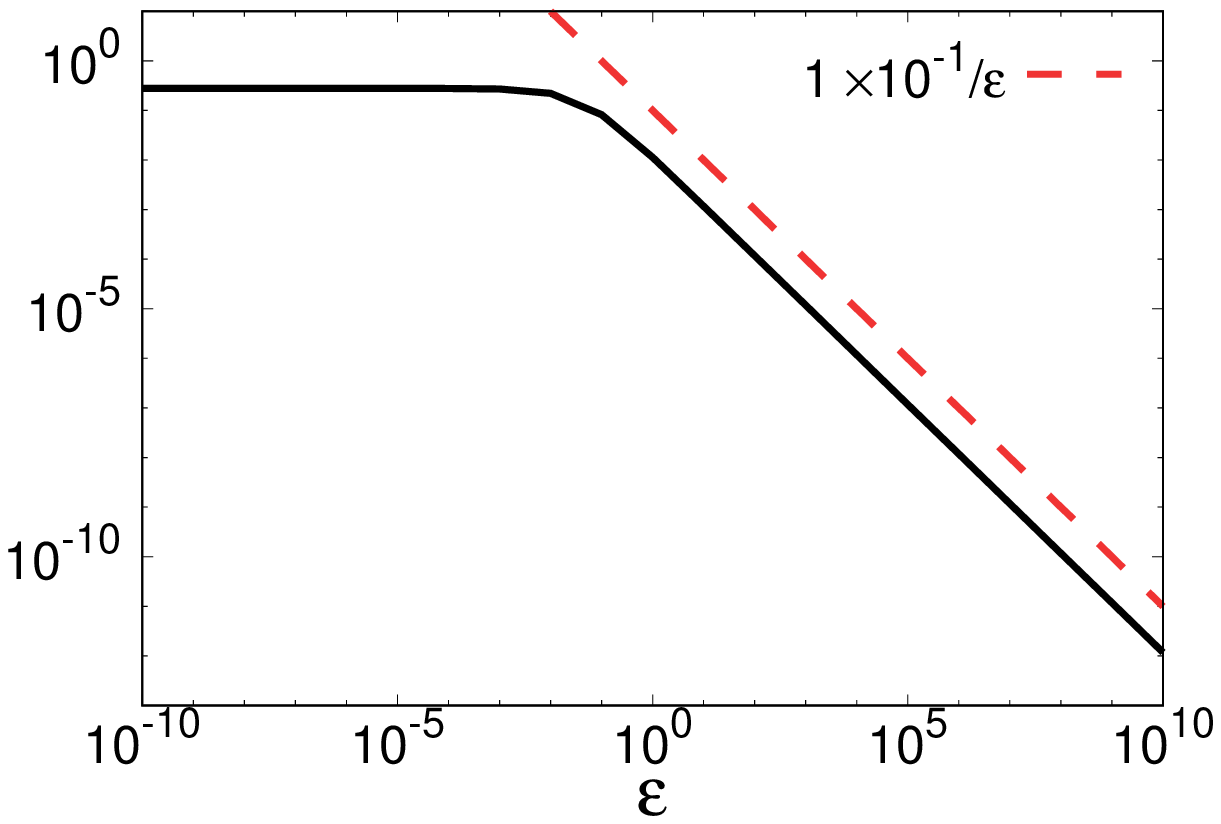}
\end{minipage}
\begin{minipage}{0.47\hsize}\centering
\includegraphics[width=\linewidth,bb=50 50 410 302]{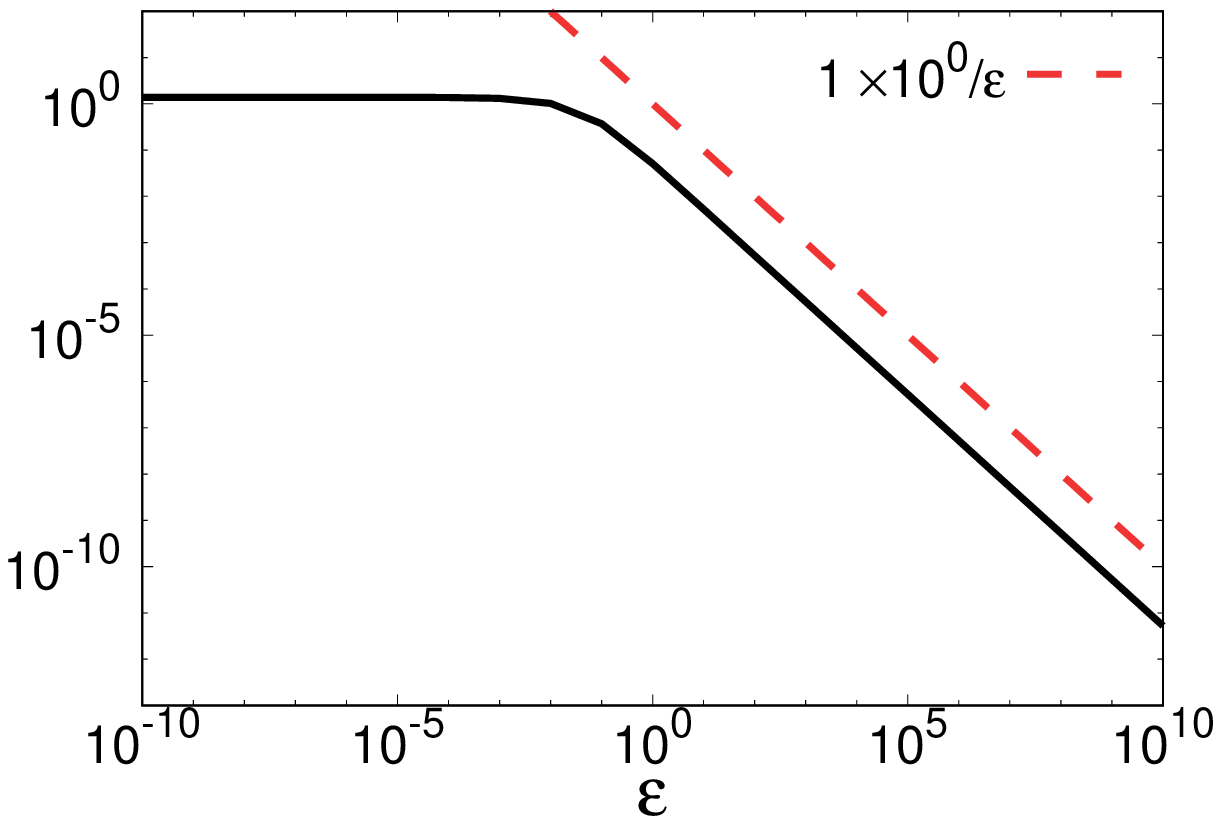}
\end{minipage}\caption{
$\|u_\eps-u_\PP\|_{\Lo^n}$ (left, solid line) and 
$\|\nabla(u_\eps-u_\PP)\|_{\Lo^{n\times n}}$ (right, solid line)
as functions of $\eps$.
}\label{epu_fig}\end{figure}

\begin{figure}[H]
\centering\begin{minipage}{0.47\hsize}\centering
\includegraphics[width=\linewidth,bb=50 50 410 302]{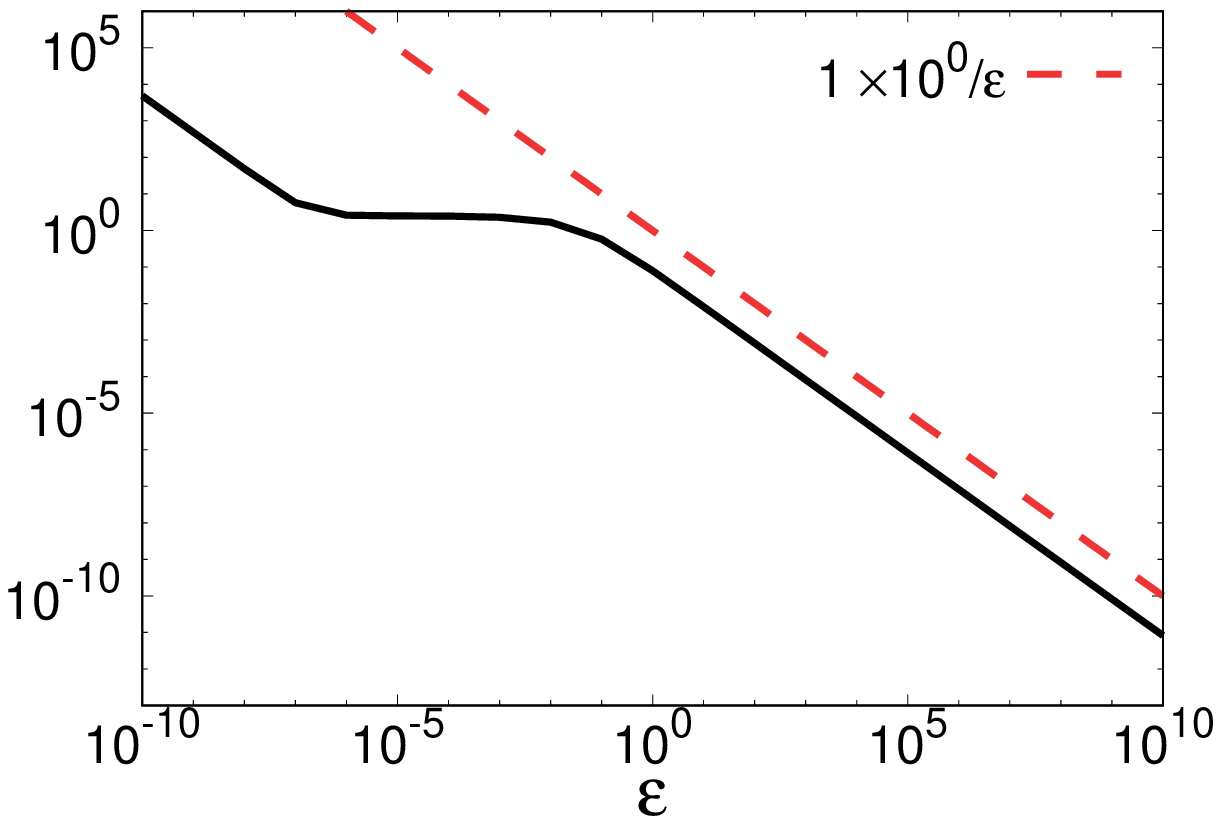}
\end{minipage}
\begin{minipage}{0.47\hsize}\centering
\includegraphics[width=\linewidth,bb=50 50 410 302]{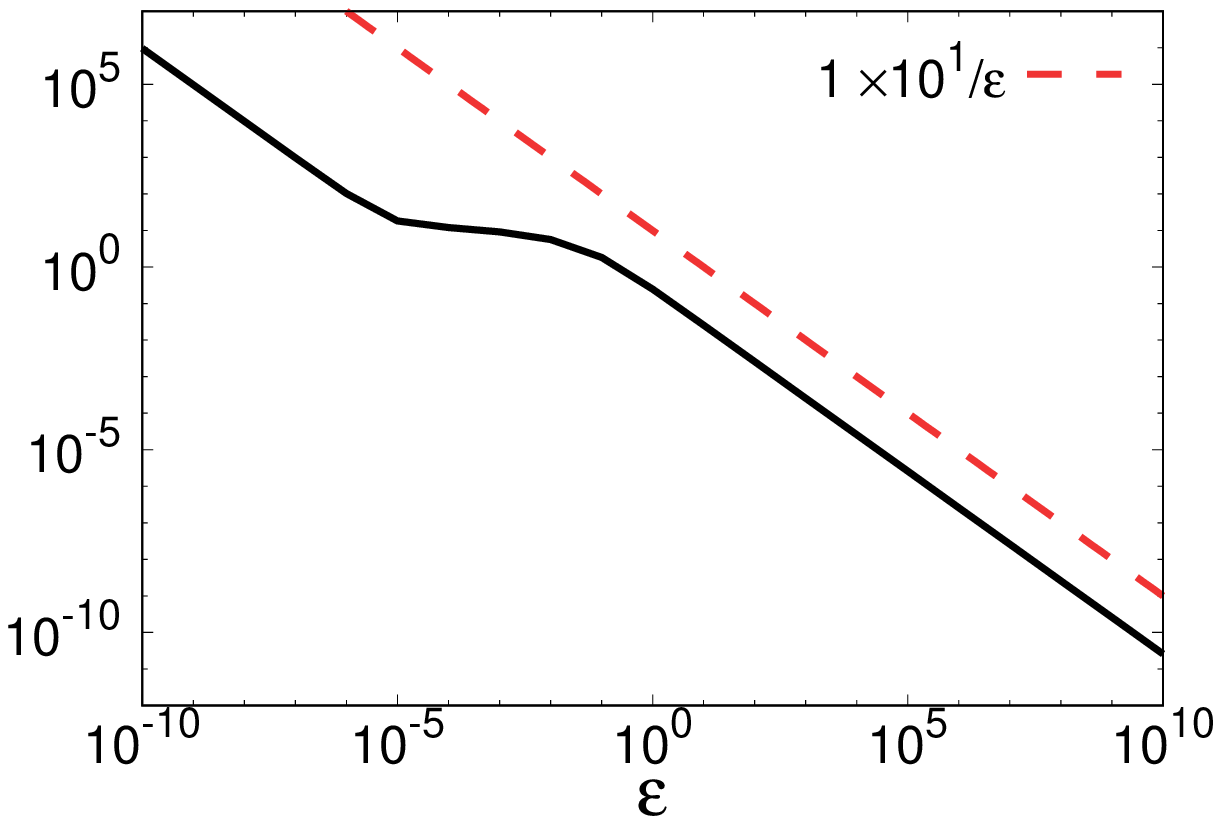}
\end{minipage}\caption{
$\|p_\eps-p_\PP\|_{\Lo}$ (left, solid line) and 
$\|\nabla(p_\eps-p_\PP)\|_{\Lo^n}$ (right, solid line)
as functions of $\eps$.}
\label{epp_fig}\end{figure}

Next we compute the error estimate between the numerical solutions of 
(ES$_1$) and (PP$_1$). The numerical errors
$\|u_\eps-u_\PP\|_{\Lo^n}$, $\|\nabla(u_\eps-u_\PP)\|_{\Lo^{n\times n}}$,
$\|p_\eps-p_\PP\|_{\Lo}$ and $\|\nabla(p_\eps-p_\PP)\|_{\Lo^n}$
are shown in Fig. \ref{epu_fig} and Fig. \ref{epp_fig}.
Based on these values, we have fitted a constant $c$ such that
$\|u_\eps-u_\PP\|_{\Hb^n}\sim c/\eps$ and
$\|p_\eps-p_\PP\|_{\Hb}\sim c/\eps$ for $\eps$ large.
Fig. \ref{epu_fig} and Fig. \ref{epp_fig} indicate that 
there exists a constant $c$ such that
$\|u_\eps-u_\PP\|_{\Hb^n}\le c/\eps$ and
$\|p_\eps-p_\PP\|_{\Hb}\le c/\eps$,
as expected from Theorem \ref{ep_conv}.

\begin{figure}[H]
\centering\begin{minipage}{0.47\hsize}\centering
\includegraphics[width=\linewidth,bb=50 50 410 302]{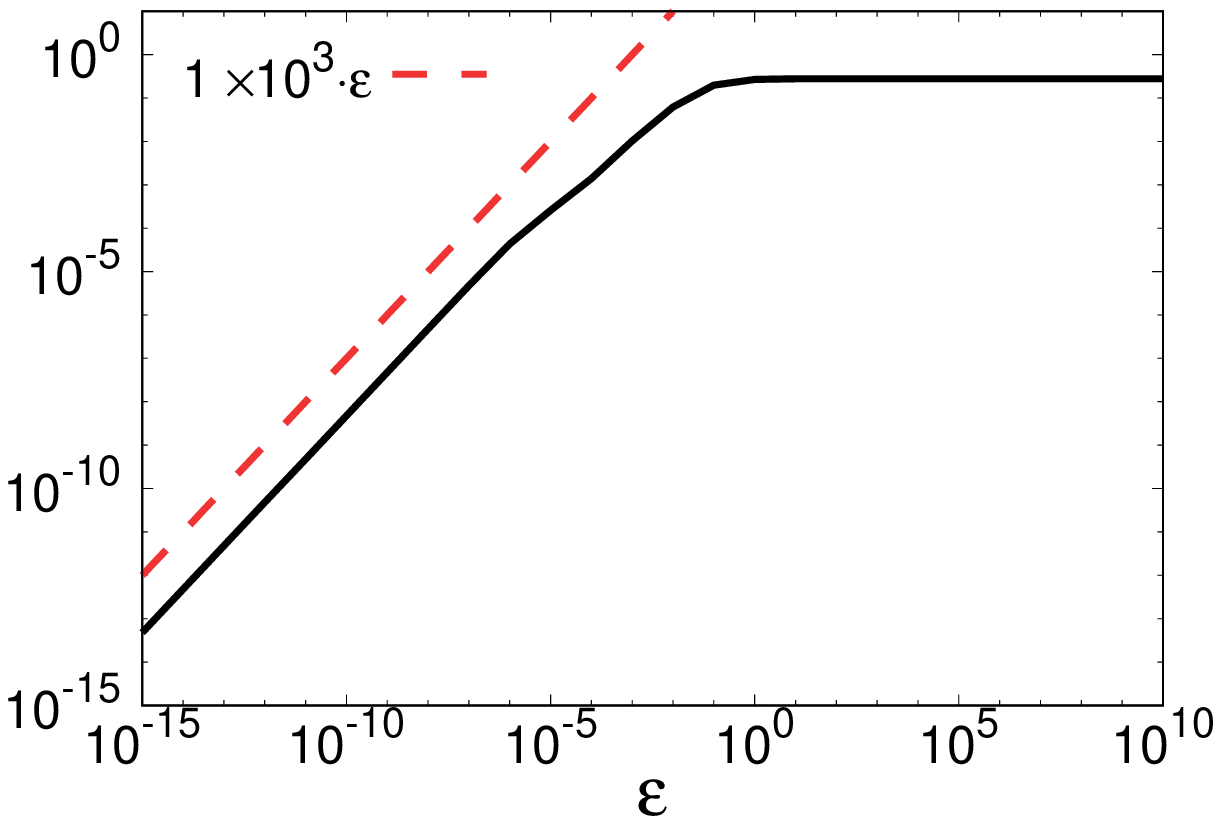}
\end{minipage}
\begin{minipage}{0.47\hsize}\centering
\includegraphics[width=\linewidth,bb=50 50 410 302]{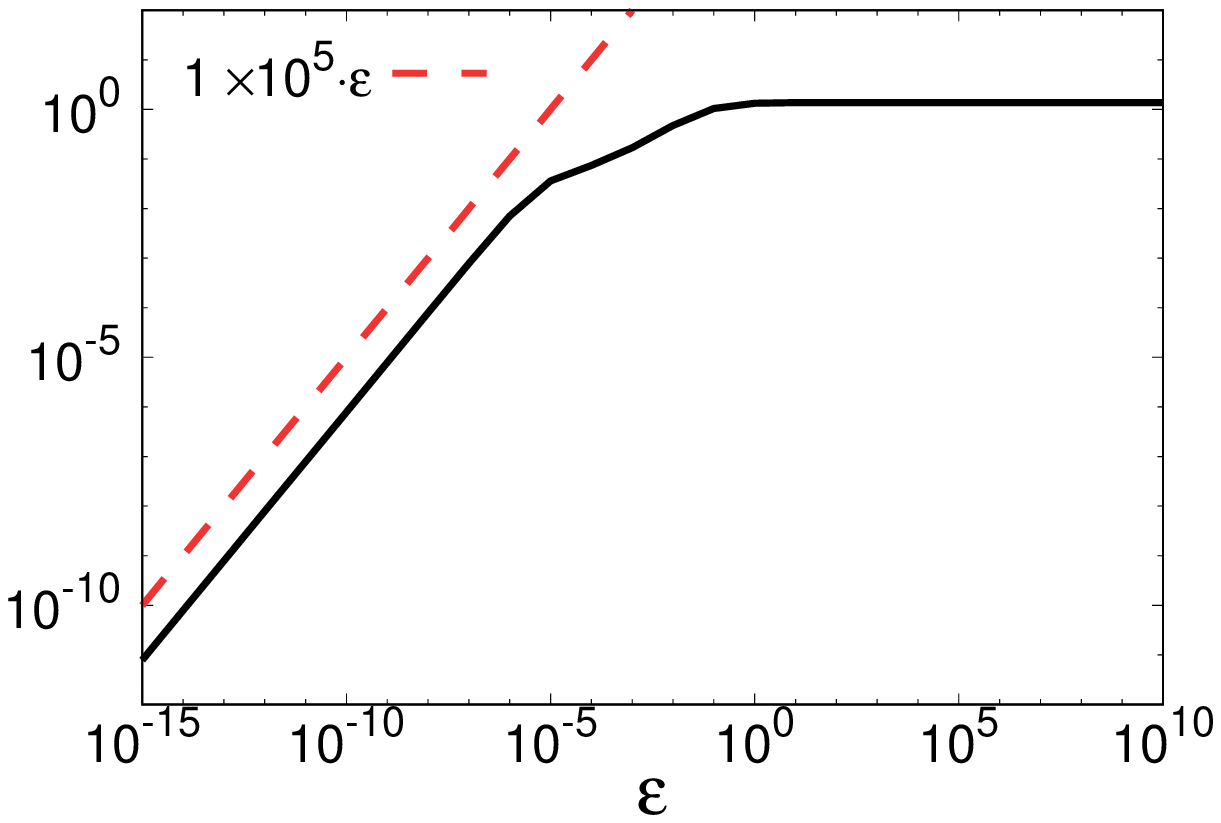}
\end{minipage}\caption{
$\|u_\eps-u_\St\|_{\Lo^n}$ (left, solid line) and
$\|\nabla(u_\eps-u_\St)\|_{\Lo^{n\times n}}$ (right, solid line)
as functions of $\eps$.}
\label{esu_fig}\end{figure}

\begin{figure}[H]
\centering\begin{minipage}{0.47\hsize}\centering
\includegraphics[width=\linewidth,bb=50 50 410 302]{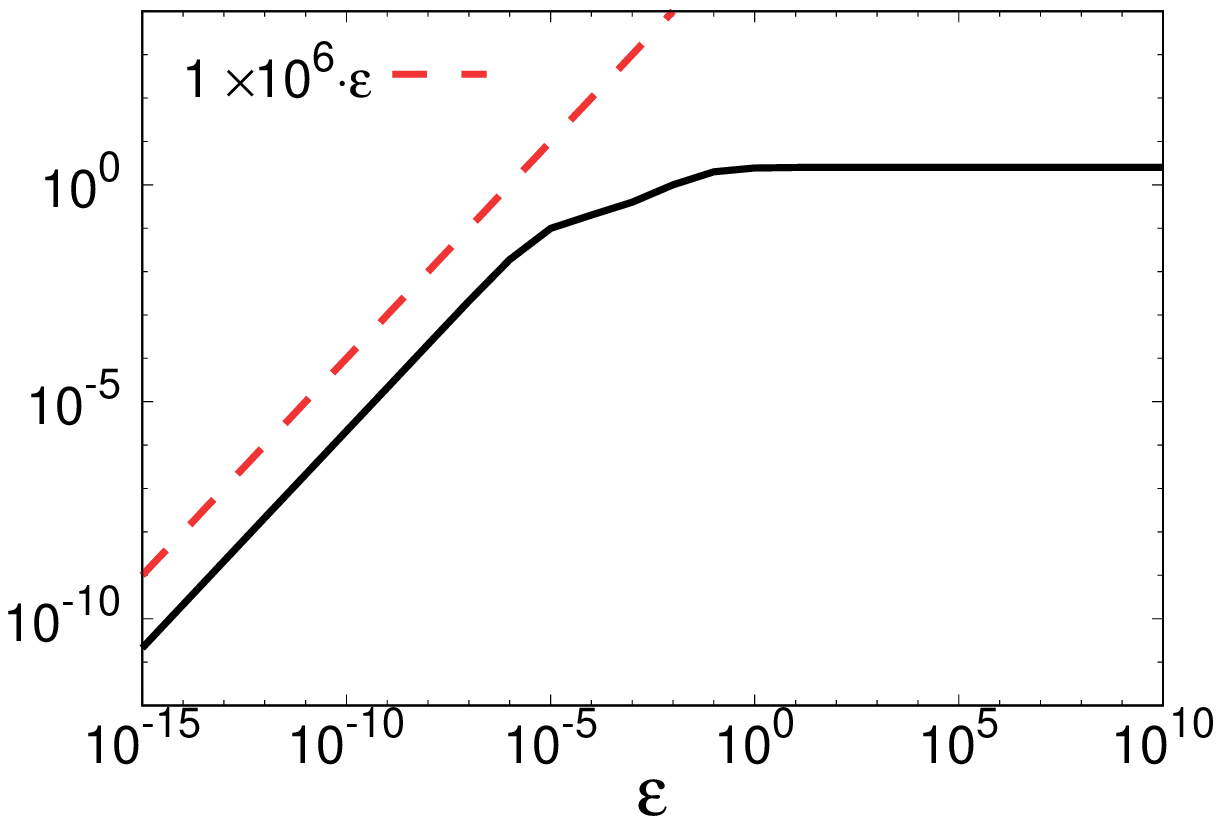}
\end{minipage}
\begin{minipage}{0.47\hsize}\centering
\includegraphics[width=\linewidth,bb=50 50 410 302]{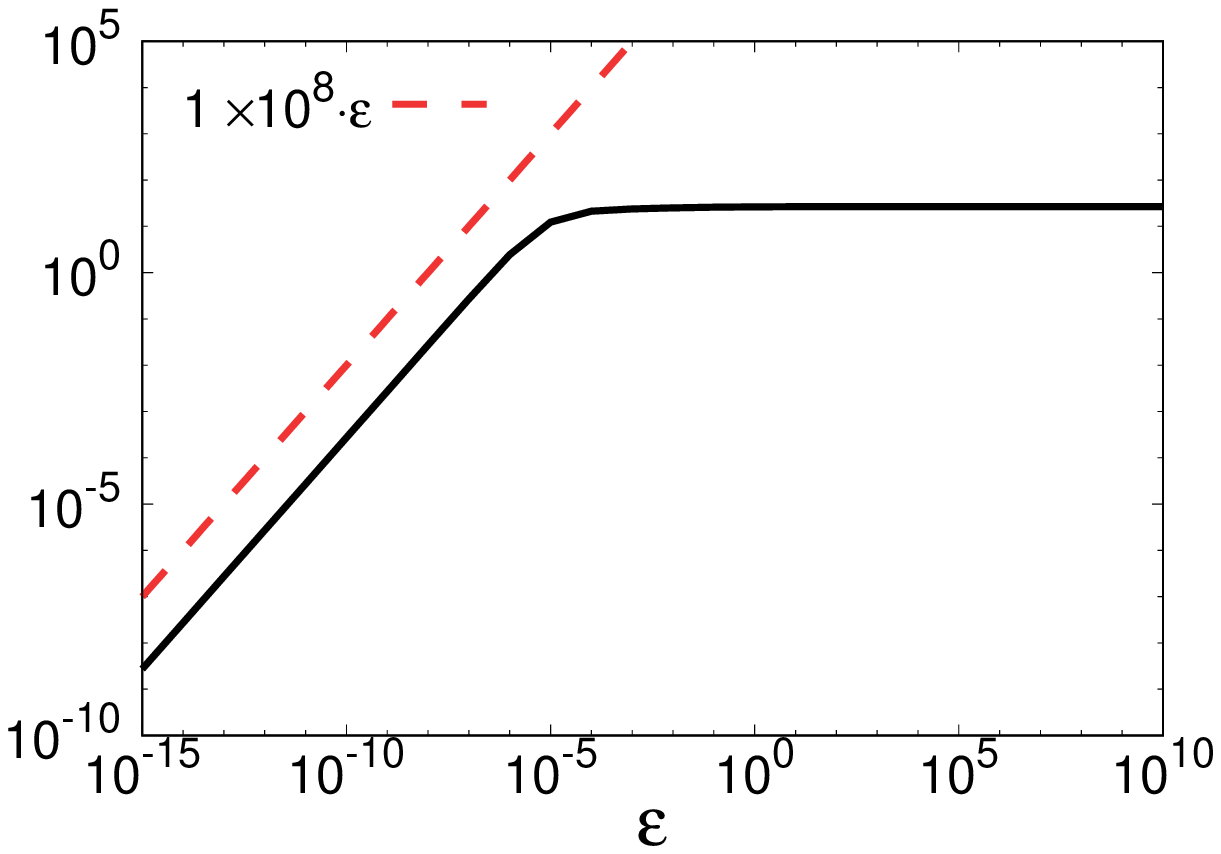}
\end{minipage}
\caption{
$\|p_\eps-p_\St\|_{\Lo}$ (left, solid line) and
$\|\nabla(p_\eps-p_\St)\|_{\Lo^n}$ (right, solid line)
as functions of $\eps$.
}\label{esp_fig}
\end{figure}

We also compute the error estimate between the problems (ES$_1$) and (S') by numerical calculation.
The numerical error estimate
$\|u_\eps-u_\St\|_{\Lo^n},~\|\nabla(u_\eps-u_\St)\|_{\Lo^{n\times n}}~,
\|p_\eps-p_\St\|_{\Lo}$ and $\|\nabla(p_\eps-p_\St)\|_{\Lo^n}$
are shown in Fig. \ref{esu_fig} and Fig. \ref{esp_fig}.
Based on these values, we have fitted a constant $c$ such that
$\|u_\eps-u_\St\|_{\Hb^n}\sim c\eps$ and
$\|p_\eps-p_\St\|_{\Lo}\sim c\eps$ for $\eps$ small.
Fig. \ref{esu_fig} and Fig. \ref{esp_fig} indicate that 
there exists a constant $\tilde{c}$ such that
$\|u_\eps-u_\St\|_{\Hb^n}\le \tilde{c}\sqrt{\eps}$ and
$\|p_\eps-p_\St\|_{\Lo}\le \tilde{c}\sqrt{\eps}$,
as expected from Theorem \ref{esn_conv}.

%\begin{acknowledgements}
%\end{acknowledgements}

\bibliographystyle{spmpsci}      % mathematics and physical sciences
\bibliography{estokes_JJIAM}   % name your BibTeX data base

%\begin{thebibliography}{}
%\end{thebibliography}

\appendix
\section*{Appendix}\label{prev_thms}

\setcounter{equation}{0}
\renewcommand{\theequation}{\Alph{section}.\arabic{equation}}
\setcounter{figure}{0}
\renewcommand{\thefigure}{\Alph{section}.\arabic{figure}}
\setcounter{table}{0}
\renewcommand{\thetable}{\Alph{section}.\arabic{table}}

Theorem \ref{estokes_thm}, \ref{bb}, \ref{es_conv},
and Proposition \ref{sp_prop} are generalizations of several theorems
stated in \cite{prev}.
In this appendix, however, we give their proofs for the readers'
convenience.
We define a continuous coercive bilinear form depending on $\eps$
and prove Theorem \ref{estokes_thm} by the Lax--Milgram Theorem.

%---------------------------------------------
%---                         Proof of estokes_thm             ---
%---------------------------------------------
\medskip
\noindent $\textit{Proof of Theorem}$ \ref{estokes_thm}.
We take arbitrary $u_1\in\Hb^n$ with $\gamma_0 u_1=u_b$.
Since $\dvg:\Ho^n\rightarrow\Lo/\R$ is surjective
\cite[Corollary 2.4, 2${}^\circ$]{Girault},
there exists $u_2\in\Ho^n$ such that
$\dvg u_2=\dvg u_1$. We put 
\begin{align}\label{def_u0}
  u_0:=u_1-u_2,
\end{align}
and note that $\gamma_0 u_0=u_b$ and $\dvg u_0=0$.
To simplify the notation, we set
$u:=u_\eps-u_0\in\Ho^n,p:=p_\eps-p_b\in Q$, and define 
$f\in\Hi^n$ and $g\in Q^*$ by 
\begin{align}\label{def_fg}\begin{array}{rll}
  \langle f,v\rangle
  &{\DS :=\into Fv-\into\nabla u_0:\nabla v-\into(\nabla p_b)\cdot v}
  &\mbox{for all }v\in\Ho^n,\\[8pt]
  \langle g,q\rangle
  &{\DS :=\langle G,q\rangle-\into\nabla p_b\cdot\nabla q}
  &\mbox{for all }q\in Q.
\end{array}\end{align}
Then, $(u_\eps,p_\eps)$ satisfies (ES') if and only if 
$(u,p)$ satisfies
\begin{align}\label{estokes_exist}
  \left\{\begin{array}{ll}
    {\DS
      \into\nabla u:\nabla \varphi +\into (\nabla p)\cdot\varphi
      =\langle f,\varphi\rangle  }
    & {\rm \mbox{for all }} \varphi\in\Ho^n,\\[8pt]
    {\DS
      \eps\into\nabla p\cdot\nabla\psi +\into(\dvg u)\psi
      =\eps\langle g,\psi\rangle }
    & {\rm \mbox{for all }} \psi \in Q.
  \end{array}\right.
\end{align}
Adding the equations in (\ref{estokes_exist}), we get
\[
  \left(\left( \begin{array}{c}u\\ p\end{array}\right) ,
  \left( \begin{array}{c}\varphi \\ \psi\end{array}\right) \right)_\eps
  :=\into\nabla u:\nabla\varphi + \eps\into \nabla p\cdot \nabla \psi
  +\into(\nabla p)\cdot\varphi + \into (\dvg u)\psi
  = \langle f,\varphi\rangle +\eps\langle g,\psi\rangle.
\]
We check that $(\cdot ,\cdot )_\eps$ is a continuous coercive bilinear form
on $\Ho^n\times Q$. The bilinearity and continuity of $(\cdot ,\cdot )_\eps$ are obvious.
The coercivity of $(\cdot ,\cdot )_\eps$ is obtained in the following way.
Take $(v,q)^T\in\Ho^n\times Q$. 
We have the following sequence of inequalities:
\begin{align*}\begin{array}{rl}
  \left(\left( \begin{array}{c}v\\ q\end{array}\right) ,
  \left( \begin{array}{c}v\\ q\end{array}\right) \right)_\eps
  =&{\DS \into\nabla v:\nabla v + \eps\into \nabla q\cdot \nabla q
    +\into v\cdot\nabla q+\into (\dvg v)q}\\[8pt]
  =&\|\nabla v\|^2_{\Lo}+\eps\|\nabla q\|^2_{\Lo}\\[8pt]
  \ge&{\DS {\rm min}\{ 1,\eps\}\left(\|\nabla v\|^2_{\Lo}+\|\nabla q\|^2_{\Lo}\right)}\\[8pt]
  \ge&{\DS c~{\rm min}\{ 1,\eps\}\left(\|v\|^2_{\Hb^n}+\|q\|^2_{\Hb}\right)}.
\end{array}\end{align*}
Summarizing, $(\cdot ,\cdot )_\eps$ is a continuous coercive bilinear form and
$\Ho^n\times Q$ is a Hilbert space.
Therefore, the conclusion of Theorem \ref{estokes_thm} follows from
the Lax--Milgram Theorem.
\qed%---------------------------------------------

Let $(u_\St,p_\St),(u_\PP,p_\PP)$ and $(u_\eps,p_\eps)$ be 
the solutions of (S'), (PP') and (ES'), respectively,
as guaranteed by Theorem \ref{stokes_t}, \ref{pp_thm} 
and \ref{estokes_thm}.
We show that the subtract $p_\St-p_\PP$ satisfies
\[
  \Delta(p_\St-p_\PP)=0
\]
in distributions sense. 
The weak harmonicity is the key ingredient to proving Proposition \ref{sp_prop}.

%---------------------------------------------
%---                         Proof of sp_prop                     ---
%---------------------------------------------
\medskip
\noindent $\textit{Proof of Proposition}$ \ref{sp_prop}.
First, we prove that there exists a constant $c>0$ independent of $\eps$
such that 
$
  \|u_\St-u_\PP\|_{\Hb^n}
  \le c\|\gamma_0 p_\St-\gamma_0 p_\PP\|_{H^{1/2}(\Gamma)},
$
and if $\gamma_0(p_\St-p_\PP)=0$, then $p_\PP=p_\St$.
Taking the divergence of the first equation of (S'), 
we obtain
\[
  \dvg F=\dvg(-\Delta u_\St+\nabla p_\St)=-\Delta (\dvg u_\St)+\Delta p_\St=\Delta p_\St.
\]
in distributions sense. 
Since $p_\St\in\Hb$ and $C^\infty_0(\Omega)$ is dense in $\Ho$,
it follows that
\[
\into\nabla p_\St\cdot\nabla\psi=-\into(\dvg F)\psi
\]
for all $\psi\in\Ho$.
Together with (S'), (PP') and $\Ho\subset Q$, we obtain
\begin{align}\label{s-p} \left\{\begin{array}{ll}{\DS
\into\nabla (u_\St-u_\PP):\nabla\varphi
=-\into (\nabla(p_\St-p_\PP))\cdot\varphi }
& {\rm \mbox{for all }} \varphi \in\Ho^n, \\[8pt]
{\DS
\into\nabla(p_\St-p_\PP)\cdot\nabla\psi=0}
&{\rm \mbox{for all }} \psi \in \Ho
\end{array}\right.\end{align}
from the assumption $\langle G,\psi\rangle=\into\nabla F\cdot\psi$.
Putting $\varphi:=u_\St-u_\PP\in\Ho^n$ in (\ref{s-p}), we get
\[\begin{array}{rl}
\|\nabla(u_\St-u_\PP)\|^2_{\Lo^{n\times n}}
&={\DS -\into (\nabla(p_\St-p_\PP))\cdot(u_\St-u_\PP)}\\[8pt]
&\le \|\nabla(p_\St-p_\PP)\|_{\Lo^n}\|u_\St-u_\PP\|_{\Lo^n}.
\end{array}\]
Hence, 
\begin{align}\label{usp_psp}
  \|u_\St-u_\PP\|_{\Hb^n}\le c_1\|\nabla(p_\St-p_\PP)\|_{\Lo^n}
\end{align}
holds.
From the second equation of (\ref{s-p}), we obtain
\[\begin{array}{rl}
  &\|\nabla (p_\St-p_\PP-\psi)\|^2_{\Lo^n}\\
  =&{\DS\|\nabla (p_\St-p_\PP)\|^2_{\Lo^n}+\|\nabla\psi\|^2_{\Lo^n}
  -2\into\nabla (p_\St-p_\PP)\cdot\nabla\psi}\\[8pt]
  =&\|\nabla (p_\St-p_\PP)\|^2_{\Lo^n}+\|\nabla\psi\|^2_{\Lo^n}\\
  \ge&\|\nabla (p_\St-p_\PP)\|^2_{\Lo^n}
\end{array}\]
for all $\psi\in\Ho$.
Thus we find
\[
\|\nabla (p_\St-p_\PP)\|_{\Lo^n}
\le
\inf_{\psi\in\Ho^n}\left(\|\nabla(p_\St-p_\PP-\psi)\|_{\Lo^n}\right).
\]
Since $\gamma_0$ is surjective and the space $\mbox{Ker}(\gamma_0)=\Ho$, 
$\Hb/\Ho$ and $H^{1/2}(\Gamma)$ are isomorphic,
there exists a constant $c_2>0$ such that
$\|q\|_{\Hb/\Ho}\le c_2\|\gamma_0 q\|_{H^{1/2}(\Gamma)}$ for all $q\in\Hb$.
Hence, we obtain
\begin{align}\label{psp_gpsp}\begin{array}{rl}
  \|\nabla (p_\St-p_\PP)\|_{\Lo^n}
  &\le{\DS\inf_{\psi\in\Ho^n}\|\nabla(p_\St-p_\PP-\psi)\|_{\Lo^n}}\\[8pt]
  &\le{\DS\inf_{\psi\in\Ho^n}\|p_\St-p_\PP-\psi\|_\Hb}\\[8pt]
  &=\|p_\St-p_\PP\|_{\Hb/\Ho}\\
  &\le c_2\|\gamma_0 p_\St-\gamma_0 p_\PP\|_{H^{1/2}(\Gamma)}.  
\end{array}\end{align}
Together with (\ref{usp_psp}), we obtain
$
\|u_\St-u_\PP\|_{\Hb^n}
\le c_1c_2\|\gamma_0 p_\St-\gamma_0 p_\PP\|_{H^{1/2}(\Gamma)}.
$
Moreover, if $\gamma_0(p_\St-p_\PP)=0$, then $p_\PP=p_\St$.

Next, we prove that there exists a constant $c>0$ independent of $\eps$
such that 
$
  \|u_\St-u_\eps\|_{\Hb^n}
  \le c\|\gamma_0 p_\St-\gamma_0 p_\eps\|_{H^{1/2}(\Gamma)},
$
and if $\gamma_0(p_\St-p_\PP)=0$, then $p_\PP=p_\eps$.
Let $w_\eps:=u_\St-u_\eps\in\Ho^n$ and $r_\eps:=p_\PP-p_\eps\in Q$.
By (S'), (PP') and (ES'), we obtain
\begin{align}\label{s-p2} \left\{\begin{array}{ll}{\DS
    \into\nabla w_\eps:\nabla\varphi
    +\into(\nabla r_\eps)\cdot\varphi
    =-\into(\nabla(p_\St-p_\PP))\cdot\varphi }
  & {\rm \mbox{for all }} \varphi \in\Ho^n, \\[8pt]
  {\DS
    \eps\into\nabla r_\eps\cdot\nabla\psi+\into(\dvg w_\eps)\psi=0}
  &{\rm \mbox{for all }} \psi \in Q.
  \end{array}\right.\end{align}
Putting $\varphi:=w_\eps$ and $\psi:=r_\eps$ 
and adding the two equations of (\ref{s-p2}), we get
\begin{align}\label{wr_pspw}
  \|\nabla w_\eps\|^2_{\Lo^{n\times n}}+\eps\|\nabla r_\eps\|^2_{\Lo^n}
  \le \|\nabla(p_\St-p_\PP)\|_{\Lo^n}\|w_\eps\|_{\Lo^n}
\end{align}
from $\into(\nabla r_\eps)\cdot w_\eps=-\into(\dvg w_\eps)r_\eps$.
Thus we find 
\[
  \|w_\eps\|_{\Hb^n}\le c_3\|\nabla(p_\St-p_\PP)\|_{\Lo^n}.
\]
Together with (\ref{psp_gpsp}), we obtain
\[
  \|u_\St-u_\eps\|_{\Hb^n}
  =\|w_\eps\|_{\Hb^n}
  \le c_2c_3\|\gamma_0 p_\St-\gamma_0 p_\PP\|_{H^{1/2}(\Gamma)}.
\]
Moreover, by (\ref{wr_pspw}), we obtain 
\[
  \eps\| p_\PP-p_\eps\|^2_\Lo
  =\eps\| r_\eps\|^2_\Lo
  \le c_4\|\nabla(p_\St-p_\PP)\|_{\Lo^n}\|w_\eps\|_{\Lo^n}.
\]
Hence, if $\gamma_0(p_\St-p_\PP)=0$, then $p_\PP=p_\eps$.
\qed%---------------------------------------------
\medskip

We show that the sequence $((u_\eps,p_\eps))_{\eps>0}$ is bounded 
in $\Hb^n\times(\Lo/\R)$.
By the reflexivity of $\Hb^n\times(\Lo/\R)$,
the sequence $((u_\eps,p_\eps))_{\eps>0}$ has a subsequence
converging weakly to somewhere in $\Hb^n\times(\Lo/\R)$.
It is sufficient to check that the limit satisfies (S').
Since the solution of (S') is unique,
the sequence $((u_\eps,p_\eps))_{\eps>0}$ converges weakly.

%---------------------------------------------
%---                         Proof of bb                       ---
%---------------------------------------------
\medskip
\noindent $\textit{Proof of Theorem}$ \ref{bb}.
We take $u_b\in\Hb^n,f\in\Hi^n$ and $g\in Q^*$
as (\ref{def_u0}) and (\ref{def_fg}) in the proof of Theorem \ref{estokes_thm}. 
We put $\tilde u_\eps:=u_\eps-u_b\in\Ho^n,\tilde p_\eps:=p_\eps-p_b\in Q$.
Then we obtain
\begin{align}\label{es_hform}
  \left\{\begin{array}{ll}
    {\DS
      \into\nabla\tilde u_\eps:\nabla \varphi +\into (\nabla\tilde p_\eps)\cdot\varphi
      =\langle f,\varphi\rangle  }
    & {\rm \mbox{for all }} \varphi\in\Ho^n,\\[8pt]
    {\DS
      \eps\into\nabla\tilde p_\eps\cdot\nabla\psi +\into(\dvg \tilde u_\eps)\psi
      =\eps\langle g,\psi\rangle }
    & {\rm \mbox{for all }} \psi \in Q.
  \end{array}\right.
\end{align}
Putting $\varphi:=\tilde u_\eps,\psi:=\tilde p_\eps$ and adding the two equations of
(\ref{es_hform}), we get
\[
\|\nabla\tilde u_\eps\|^2_{\Lo^{n\times n}}+\eps\|\nabla\tilde p_\eps\|^2_{\Lo^n}
\le \|f\|_{\Hi^n}\|\nabla\tilde u_\eps\|_{\Lo^{n\times n}}
+\eps \|g\|_{Q^*}\|\nabla\tilde p_\eps\|_{\Lo^n}
\]
since $\into(\nabla\tilde p_\eps)\cdot\tilde u_\eps=-\into(\dvg\tilde u_\eps)\tilde p_\eps$.
Hence, $(\|\tilde u_\eps\|_{\Hb^n})_{0<\eps<1}$ and
$(\|\sqrt{\eps}\tilde p_\eps\|_\Hb)_{0<\eps<1}$ are bounded. 
Moreover, by Lemma \ref{pleu}, we obtain
\[
\|\tilde p_\eps\|_{\Lo/\R}\le{c}(\|\nabla\tilde u_\eps\|_{\Lo^{n\times n}}+\|f\|_{\Hi^n}),
\]
i.e., $(\|\tilde p_\eps\|_{\Lo/\R})_{0<\eps<1}$ is bounded.
By Theorem \ref{ep_conv}, $(\|u_\eps\|_{\Hb^n})_{\eps\ge 1}$ and
$(\|\tilde p_\eps\|_{\Lo/\R})_{\eps\ge 1}$ are bounded, and thus
$(\|u_\eps\|_{\Hb^n})_{\eps>0}$ and $(\|\tilde p_\eps\|_{\Lo/\R})_{\eps>0}$
are bounded.

Since $\Hb^n\times(\Lo/\R)$ is reflexive and
$(\tilde u_\eps,[\tilde p_\eps])_{0<\eps<1}$ is bounded in $\Hb^n\times(\Lo/\R)$,
there exist $(u,p)\in\Hb^n\times(\Lo/\R)$ and a subsequence of pairs 
$(\tilde u_{\eps_k}, \tilde p_{\eps_k})_{k\in\N}\subset \Ho^n\times Q$
such that
\[
\tilde u_{\eps_k}\rightharpoonup u~{\rm weakly~\mbox{in }}\Hb^n,~
[\tilde p_{\eps_k}]\rightharpoonup p~{\rm weakly~\mbox{in }}\Lo/\R\quad{\rm as~}k\rightarrow\infty.
\]
Hence, from (\ref{es_hform}) with $\eps:=\eps_k$, taking $k\rightarrow\infty$,
we obtain
\begin{align}\label{eslimit}
\left\{ \begin{array}{ll}{\displaystyle
\into\nabla u:\nabla\varphi +\langle \nabla p,\varphi\rangle
=\langle f,\varphi\rangle }
& {\rm \mbox{for all }} \varphi\in\Ho^n\\[8pt]
{\displaystyle \into(\dvg u)\psi=0}& {\rm \mbox{for all }} \psi \in Q,
\end{array}\right.
\end{align}
where we have used that
\[
|\eps_k\into\nabla\tilde p_{\eps_k}\cdot\nabla\psi|
\le \sqrt{\eps_k}\|\sqrt{\eps}\tilde p_\eps\|_\Hb\|\psi\|_\Hb\rightarrow 0,
\]
\[
\into\nabla\tilde p_{\eps_k}\cdot\varphi
=-\into[\tilde p_{\eps_k}]\dvg\varphi
\rightarrow-\into p\dvg\varphi=\langle\nabla p,\varphi\rangle
\]
as $k\rightarrow\infty$.
By (\ref{def_fg}), the first equation of (\ref{eslimit}) implies that
\[
\into\nabla (u+u_b):\nabla\varphi +\langle \nabla (p+p_b),\varphi\rangle
=\into F\cdot\varphi
\]
for all $\varphi\in\Ho^n$.
From the second equation of (\ref{eslimit}) and 
$C^\infty_0(\Omega)\subset Q$, $\dvg(u+u_b)=0$ follows. 
Hence, we obtain that $(u+u_b,p+[p_b])$ satisfies (S'), i.e.,
$u_\St=u+u_b$ and $p_\St=p+[p_b]$. Then we have
\[
u_{\eps_k}-u_\St=u_{\eps_k}-u-u_b=\tilde u_{\eps_k}-u_\St\rightharpoonup 0
\mbox{ weakly in }\Hb^n,
\]
\[
[p_{\eps_k}]-p_\St=[p_{\eps_k}-p-p_b]=[\tilde p_{\eps_k}]-p\rightharpoonup 0
\mbox{ weakly in }\Lo/\R
\]
as $k\rightarrow\infty$.
Since any arbitrarily chosen subsequence of $((u_\eps,[p_\eps]))_{0<\eps<1}$
has a subsequence which converges to $(u_\St,p_\St)$,
we can conclude the proof.
\qed%---------------------------------------------
\medskip

Using Theorem \ref{bb} and the Rellich-Kondrachov Theorem,
it is easy to prove Theorem \ref{es_conv}.

%---------------------------------------------
%---                         Proof of es_conv                  ---
%---------------------------------------------
\medskip
\noindent $\textit{Proof of Theorem }$\ref{es_conv}.
We have from (ES') and (S') that
\[
  \left\{
  \begin{array}{ll}
    {\displaystyle
  \into\nabla (u_\eps-u_\St):\nabla\varphi +\into(\nabla(p_\eps-p_\St))\cdot\varphi =0  }
    & {\rm \mbox{for all }} \varphi\in\Ho^n,\\[8pt]
    {\displaystyle
    \eps\into\nabla p_\eps\cdot\nabla\psi+\into(\dvg u_\eps)\psi=\eps\langle G,\psi\rangle}
    & {\rm \mbox{for all }} \psi \in Q.
  \end{array}
  \right.
\]
Putting $\varphi:=u_\eps-u_\St\in\Ho^n$ and $\tilde p_\St:=p_\St-p_b\in\Hb$, 
we get
\[\begin{array}{rl}
  \|\nabla(u_\eps-u_\St)\|^2_{\Lo^{n\times n}}
  &={\DS-\into(\nabla(p_\eps-p_\St))\cdot(u_\eps-u_\St)}\\[8pt]
  &={\DS-\into(\nabla(p_\eps-p_b))\cdot(u_\eps-u_\St)
    +\into(\nabla(p_\St-p_b))\cdot(u_\eps-u_\St)}\\[8pt]
  &={\DS\into(p_\eps-p_b)\dvg(u_\eps-u_\St)
    +\into(\nabla\tilde p_\St)\cdot(u_\eps-u_\St)}\\[8pt]
  &={\DS\into(p_\eps-p_b)\dvg u_\eps
    +\into(\nabla\tilde p_\St)\cdot(u_\eps-u_\St),}
\end{array}\]
since 
$-\into(\nabla(p_\eps-p_b))\cdot(u_\eps-u_\St)=\into(p_\eps-p_b)\dvg(u_\eps-u_\St)$
and $\dvg u_\St=0$. Thus,
\begin{align}\label{ues_part}
  \|\nabla(u_\eps-u_\St)\|^2_{\Lo^{n\times n}}
  =\into(p_\eps-p_b)\dvg u_\eps+\into(\nabla\tilde p_\St)\cdot(u_\eps-u_\St).
\end{align}
Putting $\psi:=p_\eps-p_b\in Q$, we have
\[
  \eps\into\nabla p_\eps\cdot\nabla(p_\eps-p_b)+\into(\dvg u_\eps)(p_\eps-p_b)
  =\eps\langle G,p_\eps-p_b\rangle.
\]
Hence,
\begin{align}\label{pes_part}
  \eps\|\nabla (p_\eps-p_b)\|^2_{\Lo^n}
  =-\eps\into \nabla (p_\eps-p_b)\cdot\nabla p_b
  -\into(p_\eps-p_b)\dvg u_\eps+\eps\langle G,p_\eps-p_b\rangle.
\end{align}
Together with (\ref{ues_part}) and (\ref{pes_part}), we obtain
\[\begin{array}{ll}
&\|\nabla(u_\eps-u_\St)\|^2_{\Lo^{n\times n}}+\eps\|\nabla (p_\eps-p_b)\|^2_{\Lo^n}\\[8pt]
=&{\DS \into\nabla\tilde p_\St\cdot (u_\eps-u_\St)
-\eps\into \nabla (p_\eps-p_b)\cdot\nabla p_b
+\eps\langle G,p_\eps-p_b\rangle}\\[8pt]
\le &\|\nabla\tilde p_\St\|_{\Lo^n}\|u_\eps-u_\St\|_{\Lo^n}
+\eps(\|\nabla p_b\|_{\Lo^n}+\|G\|_{Q^*})\|\nabla (p_\eps-p_b)\|_{\Lo^n}.\\
\end{array}\]
By Theorem \ref{bb} and the Rellich-Kondrachov Theorem,
there exists a sequence $(\eps_k)_{k\in\N}\subset \R$
such that
\[
  u_{\eps_k}\rightarrow u_\St\mbox{ strongly in }\Lo^n
  \mbox{ as }k\rightarrow\infty.
\]
Therefore,
\[\begin{array}{ll}
  &\|\nabla(u_{\eps_k}-u_\St)\|^2_{\Lo^{n\times n}}\\
  \le&\|\nabla\tilde p_\St\|_{\Lo^n}\|u_{\eps_k}-u_\St\|_{\Lo^n}
  +\eps_k(\|\nabla\tilde p_b\|_{\Lo^n}+\|G\|_{Q^*})\|\nabla (p_{\eps_k}-p_\St)\|_{\Lo^n}\\
  \rightarrow &0
\end{array}\]
as $k\rightarrow\infty$. This implies that
\[
\|[p_{\eps_k}]-p_\St\|_{\Lo/\R}=\|p_{\eps_k}-p_\St\|_{\Lo/\R}\le c\|\nabla(u_{\eps_k}-u_\St)\|_{\Lo^{n\times n}}
\rightarrow 0\mbox{ as }k\rightarrow\infty
\]
by Lemma \ref{pleu}.
Since any arbitrarily chosen subsequence of $((u_\eps,[p_\eps]))_{0<\eps<1}$
has a subsequence which converges to $(u_\St,p_\St)$,
we can conclude the proof.
\qed%---------------------------------------------

\end{document}